\documentstyle[12pt]{article}
\newcommand{\en}{\enspace}                         
\newcommand{\bi}{\bigskip}                            
\newcommand{\me}{\medskip}                            
\newcommand{\no}{\noindent}                        
\newcommand{\be}{\begin{equation}}                   
\newcommand{\ee}{\end{equation}}                   
\newcommand{\bea}{\begin{eqnarray}}                                              
\newcommand{\eea}{\end{eqnarray}}                  
\newcommand{\timi}{\rm > \hspace{-1mm} \lhd} 
            
\newcommand{\R}{{\rm I} \! {\rm R}}

\begin{document}

\no                                                  
{\LARGE \bf Lie Groups of Fourier Integral Operators on Open Manifolds} 
                                                   
\bi                       
                             
\no                                                   
{\sc J\"urgen Eichhorn and Rudolf Schmid}         
\footnote{Research supported by
the Emory--Greifswald Exchange Program}

\bi                       
\bi                       
\bi
\bi                  
\bi

\no                   
{\bf Abstract:} We endow the group of invertible Fourier integral 
operators on an {\em open} manifold with the
structure of an ILH Lie group. 
This is done by establishing such  structures for
the groups of invertible pseudodifferential operators and contact transformations 
on an open manifold of bounded geometry, and gluing those together via a local section.

\bi
\no                      

\no

\bi
\no
{\bf Contents:}

1. Introduction

2. Bounded Geometry and Sobolev Diffeomorphism Groups

3. The Group of Contact Transformations

4. Contact Transformations of the Restricted Cotangent Bundle

5. Pseudodifferential and Fourier Integral Operators on Open
Manifolds

6. $({\cal U} \Psi^{0,k})_*$ as ILH Lie group
                  
7. An ILH Lie group structure for invertible Fourier integral
operators

\section{Introduction}    
                      
For finite dimensional Hamiltonian systems, the embedding as a 
coadjoint orbit is a very well known and convenient method for 
integration, at least for the construction of integrals. The same is
valid in the infinite dimensional case. But, as in the 
finite dimensional case, the main problem is to find appropriate Lie groups                   
such that the Hamiltonian system can be embedded as a coadjoint
orbit. There are not many such candidates. For example, 
considerations of completed diffeomorphism groups leads to the 
following complications. First, they have a good Hilbert manifold
structure but left multiplication and forming the inverse are only
continuous operations, i. e. they don't have a Lie group structure. Secondly,
considering the tangent space at the identity as a candidate for
a Lie algebra causes additional difficulties. Namely, forming the
Lie bracket decreases the Sobolev order, i. e. it is not a closed
operation. One way out of this difficulty has been indicated by
Omori \cite{Om}, forming the inverse limit of such groups and algebras, 
labeled by the Sobolev index. In the compact case, this is an old
and very well known story. In the open case, there arise tremendous
difficulties which have been essentially overcome e. g. in \cite{Ei3}, \cite{ES}.
There are not many other candidates for infinite dimensional Lie
groups. Another very important class are the invertible 
pseudodifferential operators $\Psi$DOs and Fourier integral operators 
FIOs on a manifold.
In the compact case, they have been established in \cite{ARS1}, \cite{ARS2} and have been applied by Adler in \cite{Ad} and Adams, Ratiu and Schmid \cite{ARS3}
to the complete integrability of the periodic 
KdV equation on the circle. In the open case, there has not yet been any approach
until today and we attack this problem in this paper.

Roughly speaking, we have an exact sequence
\[ I \rightarrow (\Psi D O)_* \rightarrow (F I O)_* \rightarrow
   {\cal D}_\theta (T^*M \setminus O) \rightarrow e , \]
where $I$ is the identity operator, $e$ is the identity diffeomorphism $( ~)_*$ denotes the invertible elements and ${\cal D}_\theta$
is the group of contact transformations. The main task is to establish
an IHL Lie group structure on the middle term. For this one establishes
such a structure for the boundary terms and carries it over to the middle
term via a local section, i. e. we have to perform 4 steps. Namely
to establish appropriate manifold group structures to the boundary
terms, to construct local sections and to establish such a structure
for the middle term. We follow the ideas in \cite{ARS1}, \cite{ARS2} (Adams, Ratiu,
Schmid).

We can do these steps, assuming bounded geometry and uniform boundedness in all
situations. We assume the manifold $M^n$ to be endowed with a
metric $g$ of bounded geometry of infinite order and such that 
\[ \inf \sigma_e (\bigtriangleup _1(g_S) 
   |_{ ( ker \bigtriangleup_1 )^\perp} ) > 0 , \]
where $\bigtriangleup_1(g_S)$ denotes the Laplace operator acting on
1--forms and $\sigma_e$ the essential spectrum, $g_S$ the Sasaki metric
on the cosphere bundle.
All constructions 
depend on $g$, but as we will show, they really only depend on the 
connected component
$comp(g)$ of $g$ in the space of all metrics of bounded geometry.

The paper is organized as follows. In section 2 we recall the 
main important facts concerning Sobolev spaces on open manifolds and
completed diffeomorphism groups. Section 3 is devoted to the proof
that for sufficiently bounded contact forms associated to a metric of 
bounded geometry satisfying the spectral assumption above the identity
component ${\cal D}^{r+1}_{\theta,0}$ of the completed group of contact
transformations is a Hilbert manifold and topological group. If 
in addition the metric $g$ satisfies the condition $(B_\infty)$ then
${\cal D}^\infty_{\theta,0}=\lim\limits_{\leftarrow}
{\cal D}^r_{\theta,0}$ is an ILH Lie group. In section 4 we apply these
constructions to the restricted cotangent bundle $T^*M \setminus 0$
and the cosphere bundle $S(T^*M)$. Then the canonical 1--form $\theta$
on $S(T^*M)$ is a contact form associated to the Sasaki metric
$g_S$ and ${\cal D}^{r+1}_{\theta,0}(S)$ is well defined. 
We
define ${\cal D}^{r+1}_{\theta,0}(T^*M \setminus 0)$ by homogeneous of
degree one extension. At the end we give several isomorphic descriptions
of the tangent space $T_{id}{\cal D}^{r+1}_{\theta,0}(T^*M \setminus 0)$.
Section 5 is devoted to the general notions and theorems concerning
uniform $\Psi$DO's and FIO's on open manifolds of bounded geometry. We
restrict ourself to FIO's of order $q$ , $-\infty \leq q \leq +\infty$,
whose homogeneous canonical relations are graphs
$\Gamma(f)$, where $f \in {\cal D}^{r+1}_{\theta,0}(S)$ and denote these sets by
${\cal U} F^q (f)$ . Let
${\cal U} F^q (r+1) = \bigcup\limits_{f \in {\cal D}^{r+1}_{\theta,0}}
{\cal U} F^q (f)$, and 
${\cal U} F^{q,k} (r+1) = {\cal U} F^q (r+1)/{\cal U} F^{-k-1} (r+1) $.
Similarly we have spaces ${\cal U} F^{q,k} (\infty)$ based on
${\cal D}^{r+1}_{\theta,0}(T^*M \setminus 0)$. Denote by $( ~)_*$
the group of invertible elements. Then we construct a local section
$\sigma$ of the exact sequence
\[ (E S) \quad \quad I \longrightarrow  ({\cal U} \Psi^{0,k})_* 
   \begin{array}{c}
     j \\[-2ex] \longrightarrow \\[-2ex] {}
   \end{array}         
   ({\cal U} F^{0,k})_* 
   \begin{array}{c}
     \pi \\[-2ex] \longrightarrow \\[-2ex] {}
   \end{array}         
   {\cal D}^\infty_{\theta,0} \longrightarrow e . \] 
Here $j$ is the inclusion and $\pi$ the map ($\Gamma$ with canonical 
relation $\Gamma(f)$) $\mapsto f$ ). This section is given by (\ref{section}).
Using certain Sobolev uniform structures, we obtain in section 6 the Hilbert
Lie groups $({\cal U}\Psi^{0,k,S})_*$. Then 
$({\cal U}\Psi^{0,k})_* = \lim\limits_{\leftarrow} 
({\cal U}\Psi^{0,k,S})_*$ has the structure of an ILH Lie group. Hence
the outer terms in $(ES)$ have such a structure. We use these, the 
local section $\sigma$ and several group theoretic constructions to 
establish in section 7 an ILH Lie group structure for 
$({\cal U} F^{0,k})_*$. Section 7 is strongly modeled by \cite{ARS2}, but 
nevertheless the openness of the underlying manifold always requires 
additional considerations, e. g. the ILH Lie algebras are quite different
from the ones in the compact case. The main result of the paper is
theorem 7.12.

We present in this paper a big class of ILH Lie groups for the 
application of the coadjoint orbit method for the integration of infinite
dimensional Hamiltonian systems.

\section{Bounded Geometry and Sobolev Diffeomorphism Groups}

\setcounter{equation}{0}  

We give a short summary of the basic facts. For details and proofs
we refer to \cite{Ei2}, \cite{Ei3}, \cite{EF1}.
Let $(M^n,g)$ be an open Riemannian manifold. Consider the following conditions $(I)$ and
$(B_k)$.

\me 
\begin{center}
\begin{tabular}{cll}
$(I)$ & {} & $r_{inj} (M,g) = \inf_{x \in M} r_{inj}(x) > 0$ , \\ 
{} \\ 
$(B_k)$ & {} & $|\nabla^i R| \le C_i, \en 0 \le i \le k$ , 
\end{tabular}
\end{center}
\me           

\no 
where $r_{inj}$ denotes the injectivity radius and $R$ the curvature.
We say $(M^n,g)$ has bounded geometry up to order $k$ if it satisfies
$(I)$ and $(B_k)$.

\me
\no
{\bf Lemma 2.1.} {\it The condition (I) implies completeness.}
\hfill $\Box$

\me
Let $0 \le k \le \infty$ and $M^n$ be open. Then there always exists
a metric $g$ satisfying $(I)$ and $(B_k)$, i. e. there is no topological 
obstruction against metrics of bounded geometry of any order.

Let $T^u_v$ be the bundle of $u$--fold covariant and $v$--fold
contravariant tensors and define
\bea
  \Omega^p_r (T^u_v,g) &=&
  \Big\{ t \in C^\infty (T^u_v) \en \Big| \en
  |t|_{g,p,r} \en :=   \nonumber \\
  &=&  \Big( \int \sum^r_{i=0} 
  |\nabla^i t|^p_{g,x} dvol_x(g)
  \Big)^{\frac{1}{p}} < \infty \Big\} . \nonumber
\eea
Let $\overline{\Omega}^{p,r}(T^u_v,g)$ be the
completion of $\Omega^p_r(T^u_v,g)$ 
with respect to $|\en|_{g,p,r}$  ,                  
$\stackrel{\circ}{{\Omega}}{}^{p,r}(T^u_v,g)$ the completion of 
$C^\infty_c (T^u_v)$ with respect to  $|\en|_{g,p,r}$  
and $\Omega^{p,r}(T^u_v,g)$ the space of all distributional tensor fields $t$
with $|t|_{g,p,r} < \infty$. Then we have
\[ \stackrel{\circ}{{\Omega}}{}^{p,r}(T^u_v,g) \subseteq             
   \overline{\Omega}^{p,r}(T^u_v,g) \subseteq \Omega^{p,r}(T^u_v,g) . \]

\me                       
\no                       
{\bf Lemma 2.2.} {\it If g satisfies (I), $(B_k)$, then
\[ \stackrel{\circ}{{\Omega}}{}^{p,r}(T^u_v,g) =            
   \overline{\Omega}^{p,r}(T^u_v,g) = \Omega^{p,r}(T^u_v,g),
   \en 0 \le r \le k+2 . \]
}                          
                          
\me                       
\no
Define                    
\bea                      
  {}^{b,m}|t|_g & := & \sum^m_{i=0} \sup_{x \in M} |\nabla^it|_{g,x} ,
  \nonumber \\            
  {}^b_m \Omega (T^u_v,g) & := & \Big\{ t \in C^\infty (T^u_v)
  \en \Big| \en {}^{b,m}|t|_g < \infty \Big\} , \nonumber
\eea                   
${}^{b,m} \Omega (T^u_v,g)$ the completion of  ${}^b_m \Omega (T^u_v,g)$ 
with respect to ${}^{b,m}|\en|_g$ and 
${}^{b,m} \stackrel{\circ}{\Omega} (T^u_v,g)$ 
the completion of 
$C^\infty_c(T^u_v)$ with respect to ${}^{b,m}|\en|_g$. Then
${}^{b,m} \Omega (T^u_v,g) = \big\{ t \en \big| \en t \en 
\mbox{is a} \en C^m 
\en \mbox{tensor field and} \en {}^{b,m}|t|_g < \infty \big\} $.

If $(E,h,\nabla) \rightarrow (M^n,g)$ is a Riemannian vector bundle
over $(M^n,g)$ with metric connection $\nabla$ then we make the 
analogous definitions, e. g. for $s \in C^\infty_c(E)$
\[ |s|_{p,r} := \left( \int \sum^r_{i=0} |\nabla^is|^p_x
   dvol_x(g) \right)^{\frac{1}{p}} \]
and obtain analogous spaces 
 
$\stackrel{\circ}{\Omega}{}^{p,r}(E,g,h,\nabla)$ , 
$\overline{\Omega}^{p,r}(E,g,h,\nabla)$, 
$\Omega^{p,r}(E,g,h,\nabla)$, 
${}^{b,m}\stackrel{\circ}{{\Omega}}(E,g,h,\nabla)$,
${}^{b,m}\Omega(E,g,h,\nabla)$.

For $(E,h,\nabla)$  there is an analogous condition 
$(B_k(E,\nabla))$ as for $(B_k(M,g))$,

\me 
\begin{center}
\begin{tabular}{cll}
$(B_k(E,\nabla))$ & {} & $|\nabla^i R^E| \le C_i, \en 0 \le i \le k$ , 
\end{tabular}
\end{center}
\me           
where $R^E$ denotes the curvature of $(E, \nabla)$. 
The lemma 2.2 remains true
correspondingly,
\be
   \stackrel{\circ}{{\Omega}}{}^{p,r}(E,g,h,\nabla) =            
   \overline{\Omega}^{p,r}(E,g,h,\nabla) = \Omega^{p,r}(E,g,h,\nabla),
   \en 0 \le r \le k+2 . 
\ee
if $(M^n,g)$ satisfies $(I)$ and $(B_k)$.

\me
\no
{\bf Proposition 2.3.} {\it Let $(E,h,\nabla) \rightarrow (M^n,g)$
be a Riemannian vector bundle satisfying (I), $(B_k(M^n,g))$ and  
$(B_k(E,\nabla))$.

\no
{\bf a.} Assume $k \ge r, k \ge 1, r-\frac{n}{p} \ge s-\frac{n}{q},
r \ge s, q \ge p \ge 1$. Then the inclusion
\[ \Omega^{p,r}(E) \hookrightarrow \Omega^{q,s}(E) \]   
is continuous.                                         

\no
{\bf b.} If $k \ge 0, r > \frac{n}{p}+s$, then the inclusion
\[ \Omega^{p,r}(E) \hookrightarrow {}^{b,s}\Omega(E) \]   
is continuous.
}

\me
\no
We refer to \cite{EF1} for the proof. \hfill $\Box$

\me
A key role for everything below plays the module structure theorem
for Sobolev spaces.

\sloppy
\me 
\no
{\bf Theorem 2.4.}
{\it 
Let $(E_i,h_i,\nabla_i) \rightarrow (M^n,g)$
  be vector bundles with $(I)$, $(B_k(M^n,g))$, $(B_k(E_i,\nabla_i))$,
  $i=1,2$. Assume $0\le r \le r_1, r_2 \le k$. If $r=0$ assume
  \[
  \left\{ 
    \begin{array}{rcl}
      r-\frac{n}{p} & < & r_1-\frac{n}{p_1} \\
      r-\frac{n}{p} & < & r_2-\frac{n}{p_2} \\
      r-\frac{n}{p} & \le & r_1-\frac{n}{p_1} + r_2-\frac{n}{p_2} \\
      \frac{1}{p} & \le & \frac{1}{p_1} +\frac{1}{p_2}
    \end{array}
  \right\}
  \en \mbox{or} \]
\be
  \left\{ 
    \begin{array}{rcl}
      r-\frac{n}{p} &  \le  & r_1-\frac{n}{p_1} \\
      0  & <  & r_2-\frac{n}{p_2} \\
      \frac{1}{p} & \le  & \frac{1}{p_1}
    \end{array}
  \right\}
  \en \mbox{or} \en
  \left\{ 
    \begin{array}{rcl}
      0  & <  & r_1-\frac{n}{p_1} \\
      r-\frac{n}{p} &  \le  & r_2-\frac{n}{p_2} \\
      \frac{1}{p} & \le  & \frac{1}{p_2}
    \end{array}
  \right\}.
\ee 
  If $r>0$ assume $\frac{1}{p}\le \frac{1}{p_1} + \frac{1}{p_2}$ and
\be 
  \left\{ 
    \begin{array}{rcl}
      r-\frac{n}{p} & < & r_1-\frac{n}{p_1} \\
      r-\frac{n}{p} & < & r_2-\frac{n}{p_2} \\
      r-\frac{n}{p} & \le & r_1-\frac{n}{p_1} + r_2-\frac{n}{p_2}
    \end{array}
  \right\}
  \en \mbox{or} \en
  \left\{ 
    \begin{array}{rcl}
      r-\frac{n}{p} & \le & r_1-\frac{n}{p_1} \\
      r-\frac{n}{p} & \le & r_2-\frac{n}{p_2} \\
      r-\frac{n}{p} & < & r_1-\frac{n}{p_1} + r_2-\frac{n}{p_2}
    \end{array}
  \right\}.
\ee 
Then the tensor product of sections defines a continuous bilinear map
\[ \Omega^{p_1,r_1}(E_1,\nabla_1) \times \Omega^{p_2,r_2}(E_2,\nabla_2) 
   \longrightarrow
   \Omega^{p,r} (E_1 \otimes E_2, \nabla_1 \otimes \nabla_2) . \]
}   
    
\no 
We refer to \cite{EF1} for the proof. \hfill $\Box$
    
\me 
Consider now $(M^n,g), (N^{n'},h)$ open, satisfying $(I), 
(B_k)$ and $f \in C^\infty(M,N)$. Then the differential $df=f_*$ is a
section of $T^*M \otimes f^*TN$, where $f^*TN$ is endowed with the induced
connection $f^*\nabla^h$. The connections $\nabla^g$ and $f^*\nabla^h$
induce connections $\nabla$ in all tensor bundles 
$T^q_s(M) \otimes f^*T^u_v N$. Therefore $\nabla^mdf$ is well defined.
Assume $m \le k$. We denote by $C^{\infty,m}(M,N)$ the set of all
$f \in C^\infty (M,N)$ satisfying 
\[ {}^{b,m}|df| := \sum^{m-1}_{i=0} \sup_{x \in M} |\nabla^idf|_x 
   < \infty . \]
Let $Y \in \Omega(f^*TN) := C^\infty (f^*TN)$. Then $Y_x$
can be written as $(Y_{f(x)}, x)$, and we define a map
$f_Y:M \rightarrow N$ by
\[ f_Y(x) := (\exp Y)(x) := \exp Y_x := \exp_{f(x)} Y_{f(x)} . \]
Then the map $f_Y$ defines an element of $C^\infty(M,N)$. Moreover
we have

\me
\no
{\bf Proposition 2.5.} {\it
Assume $m \le k$ and 
$ {}^{b,m} |Y| = \sum\limits^m_{i=0} \sup\limits_{x \in M}
  |\nabla^iY|_x < \delta_N < r_{inj}(N)$,
$f \in C^{\infty,m}(M,N)$. Then 
\[ f_Y = \exp Y \in C^{\infty,m}(M,N) . \]
}

\me 
\no
We refer to \cite{Ei3} for the proof.  The main point is that one shows
that $|\nabla^u (d \exp Y - d(id))|$ makes sense and that
\be
  |\nabla^\mu (d \exp Y - d(id))| \le P_\mu (|\nabla^i df|, 
  |\nabla^j Y|), \en i \le \mu, \en j \le \mu+1 ,
\ee
where the $P_\mu$ are certain universal polynomials in the indicated
variables without constant terms and each term has at least one
$|\nabla^jY|, 0 \le j \le \mu+1$ as a factor. \hfill $\Box$

\me
Now we consider manifolds of maps in the $L_p$--category. Assume
that $(M^n,g), (N^{n'},h)$ are open, of bounded geometry
up to order $k, r \le m \le k, 1 \le p < \infty, r > \frac{n}{p}+1$.
Consider $f \in C^{\infty,m}(M,N)$. According to 2.3, for 
$r > \frac{n}{p}+s$
\bea 
  \Omega^{p,r}(f^*TN) & \hookrightarrow & {}^{b,s} \Omega (f^*TN) 
  \nonumber \\
  {}^{b,s}|Y| & \le & D \cdot |Y|_{p,r} , \nonumber 
\eea
where 
$ |Y|_{p,r} = \left( \int \sum\limits^r_{i=0} |\nabla^iY|^p dvol 
  \right)^\frac{1}{p} $
and $\nabla = f^* \nabla^h$ . Set for 
$\delta > 0, \delta \cdot D \le \delta_N < r_{inj}(N)/2, 1 \le p < \infty$
\bea 
   V_\delta := \Big\{ (f,g) \in C^{\infty,m} (M,N)^2 \en \Big| \en
   \mbox{there exists} \en Y \in \Omega^p_r (f^*TN) \en
   \mbox{such that} \nonumber \\
   g = f_Y = \exp Y \en \mbox{and} \en
   |Y|_{p,r} < \delta \Big\} . \nonumber
\eea

\me
\no
{\bf Theorem 2.6.} {\it Under the conditions above
${\cal V} := \{V_\delta\}_{0<\delta<r_{inj}(N)/2D}$ is a basis for a
metrizable uniform structure ${\cal A}^{p,r}(C^{\infty,m}(M,N))$.
}

\me
\no
We refer to \cite{EF1} for the rather complicated proof. 
\hfill $\Box$.

\me
Let ${}^m\Omega^{p,r}(M,N)$ be the completion of $C^{\infty,m}(M,N)$
with respect to this uniform structure. From now on we assume
$r=m$ and denote $\Omega^{p,r}(M,N) := {}^r\Omega^{p,r}(M,N)$.

\me
\no
{\bf Theorem 2.7.} {\it
Let $(M^n,g), (N^{n'},h)$ be open and of bounded geometry
of order $k$, $1 \le p < \infty, k \ge r > \frac{n}{p}+1$. Then each
component of $\Omega^{p,r}(M,N)$ is a $C^{k+1-r}$--Banach manifold,
and for $p=2$ it is a Hilbert manifold.
}

\me
\no
We refer to \cite{EF1} for the proof. \hfill $\Box$

\me
Let $(M^n,g)$ be as above. A choice of an orthonormal basis in each 
$T_xM$ implies that $|\lambda|_{min}(df)$, the minimum of the absolute
value of the eigenvalues of the Jacobian of $f$, is well defined. Set
\[ {\cal D}^{p,r} := \Big\{ f \in \Omega^{p,r}(M,N) \en \Big| \en 
   f \en \mbox{is injective, surjective and} \en  
   |\lambda|_{min}(df)>0 \Big\} . \]  
                              
\me                            
\no                           
{\bf Theorem 2.8.} {\it   
${\cal D}^{p,r}$ is open in $\Omega^{p,r}(M,N)$. In particular, each 
component is a $C^{k+1-r}$--Banach manifold, and for $p=2$ it is a
Hilbert manifold.} \hfill $\Box$
 
\me
\no     
{\bf Theorem 2.9.} {\it        
Assume $(M^n,g), k, p, r$ as above.

\no
{\bf a.} Assume $f,h \in {\cal D}^{p,r}$, 
$h \in comp(id_M) \subset {\cal D}^{p,r}$ .
Then $h \circ f \in {\cal D}^{p,r}$ and 
$h \circ f \in comp(f)$.

\no
{\bf b.} Assume $f \in comp(id_M) \subset {\cal D}^{p,r}$. Then 
$f^{-1} \in comp(id_M) \subset {\cal D}^{p,r}$.

\no
{\bf c.}
${\cal D}^{p,r}_0 := comp(id_M)$ is a metrizable topological group.
}
\hfill $\Box$

\me
\no     
{\bf Theorem 2.10 ($\alpha$--lemma).} {\it        
Assume $k \ge r > \frac{n}{p}+1, f \in {\cal D}^{p,r}$. Then the right
multiplication $\alpha_f: {\cal D}^{p,r}_0 \rightarrow {\cal D}^{p,r}$,
$\alpha_f(h) = h \circ f$, is of class $C^{k+1-r}$.
}
\hfill $\Box$

\me
\no     
{\bf Theorem 2.11 ($\omega$-lemma).} {\it          
Let $k+1-(r+s)>s$, $f \in {\cal D}^{p,r+s}_0 \subset {\cal D}^{p,r}_0 $,
$r > \frac{n}{p}+1$. Then the left multiplication 
$\omega_f: {\cal D}^{p,r} \rightarrow {\cal D}^{p,r}$,
$\omega_f(h) = f \circ h$, is of class $C^s$.
}
\hfill $\Box$

\me
We defined for $C^{\infty,m}(M,N)$ a uniform structure ${\cal A}^{p,r}$.
Consider now $C^{\infty, \infty}(M,N) = \bigcap_m C^{\infty,m}(M,N)$.
Then we have an inclusion 
$i: C^{\infty, \infty}(M,N) \hookrightarrow C^{\infty,m}(M,N)$                                       
and hence a well defined uniform structure 
${\cal A}^{\infty,p,r} = (i \times i)^{-1} {\cal A}^{p,r}$. After
completion we obtain once again the manifold
$\Omega^{\infty,p,r}(M,N)$, where $f \in \Omega^{\infty,p,r}(M,N)$
if and only if for every $\varepsilon > 0$ there exist an 
$\tilde{f} \in C^{\infty, \infty}(M,N)$ and a
$Y \in \Omega^{p,r}(\tilde{f}^*TN)$ such that $f = \exp Y$ and
$|Y|_{p,r} \le \varepsilon$. Moreover, each component of 
$\Omega^{\infty,p,r}(M,N)$ is a Banach manifold and 
$T_f \Omega^{\infty,p,r}(M,N) = \Omega^{p,r}(f^*TN)$. As above we
set
\[ {\cal D}^{\infty,p,r} = \Big\{ f \in \Omega^{\infty,p,r}(M,N) 
   \en \Big| \en 
   f \en \mbox{is injective, surjective and} \en  
   |\lambda|_{min}(df)>0 \Big\} . \]  

\me
\no
{\bf Theorem 2.12.} {\it
Assume the conditions for defining ${\cal D}^{p,r}$. Then }
\[ {\cal D}^{\infty,p,r}_0 = {\cal D}^{p,r}_0 . \] 
We refer to \cite{ES}, p. 163 for the proof. \hfill $\Box$

\section{The Group of Contact Transformations}

\setcounter{equation}{0}  

From now on we restrict ourselves to $p=2$ and write 
$ {\cal D}^r_0 \equiv {\cal D}^{2,r}_0 $. Moreover, we have to consider
$q$--forms with values in a vector bundle $E$, i. e. elements of
$\Omega^{q,2,r} (E) \equiv  \Omega^{2,r} (\Lambda^q T^* M \otimes E) $.
Sections of $E$ are simply 0--forms with values in $E$ . 
Usual forms on $M$ are forms with values in $M \times R \rightarrow M$
and we write simply $\Omega^{q,2,r} \equiv \Omega^{q,2,r}(M)  $.

In \cite{ES} we studied the group ${\cal D}^{r+1}_{\omega,0}$ of form
preserving diffeomorphisms $f \in {\cal D}^{r+1}_0$, 
$f^*\omega = \omega$, $\omega$ a symplectic or volume form. 
We proved the following

\me
\no
{\bf Theorem 3.1.} {\it
Assume $(M^n,g)$ with $(I)$, $(B_\infty)$, $\omega \in {}^{b,m} \Omega^q$
for all $m$, closed and strongly nondegenerate, $q=n$ or $q=2$,
$\inf \sigma_e (\bigtriangleup_1 |_{(ker \bigtriangleup_1)^\perp}) > 0$.
Let ${\cal D}^\infty_{\omega,0} = \lim\limits_\leftarrow 
{\cal D}^r_{\omega,0}$.
Then $\Big\{ {\cal D}^\infty_{\omega,0}, {\cal D}^r_{\omega,0} 
| r > \frac{n}{2}+1 \Big\}$ is an ILH Lie group in the sense of \cite{Om},
\cite{ARS3}  and the Lie algebra of ${\cal D}^\infty_{\omega,0}$ consists of
divergence free $(q=n)$ or locally Hamiltonian $(q=2)$ vector fields $X$,
respectively, with $|X|_{2,r} < \infty$ for all $r$. 
}
\hfill $\Box$

\me
\no
Here strongly nondegenerate means that $\inf\limits_{x \in M} |\omega|_x > 0$.

\me 
A similar theorem for the group of contact transformations would be
desirable and will be necessary for the constructions in sections 5 and 6.
As a result of our efforts, such a theorem can be established but it is 
once again a long, complicated story and will appear together with results on 
other
diffeomorphism groups in \cite{EF2}. Hence we only sketch the proof here.

Let $(M^{2n+1}, g_0,\theta)$ be an oriented Riemannian contact manifold, 
$g_0$ satisfying $(I)$ and $(B_{k+2})$ and $\theta$ a contact form. We
assume additionally ${}^{b,k+3}|\theta|_{g_0} < \infty$. $\theta$ is
a 1--form with $\mu := \theta \wedge (d \theta)^n \not= 0$ everywhere
and we assume that $\mu$ coincides with the given orientation. $\theta$
defines the Reeb vector field $\xi$ on $M$, 
$\theta (\xi) = 1$, $i_\xi d \theta =0$. Denote $E= ker \theta$.
Clearly $TM= R \xi \oplus E$. A Riemannian metric $g$ is called 
associated to $\theta$ if there exists a (1,1) tensor field $\varphi$
on $M$ such that for any vector fields $X,Y$ on $M$ we have

1) $g(X,\xi) = \theta (X)$ ,

2) $\varphi^2 = - I + \theta \otimes \xi$ ,

3) $d \theta (X,Y) = g (X, \varphi Y)$ . 

\no
These conditions imply

4) $g(\xi, \xi) =1$ ,

5) $E \perp \xi $ ,

6) $\varphi(\xi) = 0 , \varphi(E)=E$ ,

7) $d \theta (\varphi X, \varphi Y) = d \theta (X, Y)$ ,

8) $g(X,Y) = \theta (X) \cdot \theta (Y) + d \theta (\varphi X, Y)$ .

\no
Given $g_0$ with $(I)$ and $(B_{k+2})$, $\theta$ as above, we want to 
construct 
a  metric $g$ of bounded geometry associated to $\theta$. This 
can be achieved as follows: 

\me
\no
{\bf Proposition 3.1.} {\it
Assume $(M^{2n+1},g_0,\theta)$ as above. Then there exists a metric $g$
associated to $\theta$ satisfying (I), $(B_k)$ and 
${}^{b,k}|\theta|_g < \infty$ .
}

\me
\no
{\bf Proof.} We sketch the simple proof. Start with $g_0$ and define
$h$ by $h(X,Y) = g_0(-X+\theta(X)\xi, -Y+\theta(Y)\xi)+\theta(X)\theta(Y)$.
Then $\xi \perp_h ker \theta$ and $|\xi|_h=1$. Let
$(X_1, \dots, X_{2n}, \xi)$ be a local orthonormal basis with respect to
$h$ and write 
\[ ((d \eta)_{ij}) = (d \eta(X_i,Y_j)) = F \cdot G , \]
where $F$ is orthonormal and $G$ is symmetric and positive definite. Then,
according to
\cite{Bl2}, 
$ \left( \begin{array}{cc}
  G & 0 \\  
  0 & 1         
  \end{array} \right)$       
defines a Riemannian metric $g$ on $M$ and 
$ \left( \begin{array}{cc}
  F & 0 \\     
  0 & 0                  
  \end{array} \right)$       
defines a global (1,1) tensor field $\varphi$ with 
$\varphi^2 = -I + \theta \otimes \xi$ and 
$d \theta(X,Y) = g (X, \varphi Y) $.

Denote by '' $'$ '' in a symbolic notation the (euclidean) differentiation.
Then 
$ (G)' = \left( \frac{(d \eta)_{ij}}{F} \right)' =
  ((d \eta)_{ij})' \cdot \frac{1}{F} - ((d \eta)_{ij}) \cdot
  \frac{F'}{F^2}$.
Similarly for higher derivatives. This implies $(B_k)$ for $g$. Finally
the condition 
$(I)$ for $g_0$ and the fact that $g_0 \rightarrow g$ implies 
uniformly boundedness from above and below and the change of local volumes 
yields (I) for $g$. Here we use theorem 4.7 of \cite{CGT}.
\hfill $\Box$

\me
From now on we assume
$g$ with conditions $(I), (B_k)$ and properties 1) -- 8),
$k \ge r+1 $, ${}^{b,r+1}|\theta|_g < \infty$ ,
$r+1 > \frac{2n+1}{2}+2$. Consider the space
\[ {\cal F} := \Big\{ \alpha \in C^\infty (M) \en \Big| \en
   {}^{b,r+1}|\alpha| < \infty \Big\} . \]
Set for $\delta >0$
\[ V_\delta := \Big\{ (\alpha_1, \alpha_2) \in {\cal F}^2 \en \Big| \en
  |\alpha_1-\alpha_2|_{2,r+1} = \Big( \int \sum^{r+1}_{i=0} 
  |(\nabla^g)^i (\alpha_1-\alpha_2)|^2_{g,x} dvol_x(g)
  \Big)^{\frac{1}{p}} < \delta \Big\} . \]

\me
\no
{\bf Lemma 3.2.} {\it
${\cal B} = \{V_\delta\}_{\delta > 0}$ is a basis for a metrizable
uniform structure on the space ${\cal F}$}. \hfill $\Box$

\me
Let ${\cal F}^{r+1}$ be the completion of ${\cal F}$ with respect to 
${\cal B}$. Then
${\cal F}^{r+1}$ is locally contractible, hence locally arcwise connected,
hence components coincide with arc components. The elements of ${\cal F}$ are
dense in each component.

\me
\no    
{\bf Proposition 3.3.} {\it 
Let $\alpha \in {\cal F}$. Then the component of $\alpha$ is given by
\[ comp(\alpha) = \Big\{ \alpha' \in {\cal F}^{r+1} \en \Big| \en
   |\alpha-\alpha'|_{2,r+1} < \infty \Big\} = \alpha + 
   \Omega^{0,2,r+1}(M) . \]
In particular each component is open and a Hilbert manifold modeled over
$\Omega^{0,2,r+1}(M)$. 
}
\hfill $\Box$

\me
\no
{\bf Corollary 3.4.} {\it 
${\cal F}^{r+1}$ has a representation as a topological sum of its
components,
\[ {\cal F}^{r+1} = \sum_{i \in I} comp (\alpha_i) . \]
}
\hfill $\Box$

Set ${\cal F}^{r+1}_0 = \Big\{ \alpha \in comp(1) \en \Big| \en
\inf\limits_{x \in M} \alpha(x) > 0 \Big\}$. Then 
${\cal F}^{r+1}_0$ is an open
subset of $comp(1)$, in particular
$T_\alpha {\cal F}^{r+1}_0 = \Omega^{0,2,r+1}(M), 
\alpha \in {\cal F}^{r+1}_0$.
Moreover, ${\cal F}^{r+1}_0$ is a Hilbert Lie group. \\
Now we define
\[ {\cal D}^{r+1}_{\theta,0} := \Big\{ (\alpha,f) \in {\cal F}^{r+1}_0
   \timi {\cal D}^{r+1}_0 \en \Big| \en \alpha f^* \theta = \theta                 
   \Big\} , \]                                                                    
where ${\rm > \hspace{-1mm} \lhd}$
denotes the semidirect product.

\me
\no
{\bf Proposition 3.5.} {\it
${\cal D}^{r+1}_{\theta,0}$ is a closed subgroup of 
${\cal F}^{r+1}_0 \timi {\cal D}^{r+1}_0$ and a topological group.
}
\hfill $\Box$

\me
\no
{\bf Theorem 3.6.} {\it
Assume $(M^{2n+1},g,\theta)$ with (I), $(B_k)$, $g$ associated to 
$\theta$, $k \ge r+1 > \frac{2n+1}{2}+2$,
${}^{b,r+1}|\theta| < \infty$ and 
$\inf \sigma_e(\bigtriangleup_1) |_{(ker \bigtriangleup_1)^\perp} > 0$.
Then ${\cal D}^{r+1}_{\theta,0}$ is a closed $C^{k-r}$ Hilbert 
submanifold of ${\cal F}^{r+1}_0 \timi {\cal D}^{r+1}_0$.
}

\me
The sketched proof will occupy the remaining part of this section. As
usual, we will show that ${\cal D}^{r+1}_{\theta,0}$ is the preimage of a
point by a submersion.

\me
\no
{\bf Lemma 3.7.} {\it
Let $\alpha \in {\cal F}^{r+1}_0, f \in {\cal D}^{r+1}_0$. Then
\be
  \alpha f^* \theta - \theta \in \Omega^{1,2,r}
\ee
and
\be
  d \alpha \wedge f^* \theta - \alpha f^* d \theta - d \theta \in 
  \Omega^{2,2,r} .
\ee
}

\me
\no
{\bf Proof.} The proof will be based on the Lemmas 3.8 -- 3.10.

Write
\bea 
  \alpha f^* \theta - \theta & = & \alpha (f^* \theta - \theta)
  + (\alpha -1) \theta \\   
  \alpha (f^* \theta - \theta) & = & (\alpha -1) (f^* \theta - \theta)
  + f^* \theta - \theta .  
\eea 
Assume $f^* \theta - \theta \in \Omega^{1,2,r}$. We have 
$(\alpha-1) \in \Omega^{0,2,r+1}$.
The module structure theorem 2.4 applied to 
$(\alpha-1) (f^* \theta - \theta)$ gives
$(\alpha-1) (f^* \theta - \theta) \in \Omega^{1,2,r}$ and
$\alpha (f^* \theta - \theta) \in \Omega^{1,2,r}$.
Moreover, $(\alpha-1) \theta \in \Omega^{1,2,r}$ since 
$(\alpha-1) \in \Omega^{0,2,r+1}$
and ${}^{b,r+1} |\theta| < \infty$.

Hence the proof of (3.1) reduces to the following

\me
\no
{\bf Lemma 3.8.} {\it Assume $f \in {\cal D}^{r+1}_0$. Then 
$f^* \theta - \theta \in \Omega^{1,2,r}$.
}

\me 
\no
{\bf Proof.} Any $f \in {\cal D}^{r+1}_0$ has a representation 
$f=\exp X_u \circ \dots \circ \exp X_1$. We start with the simplest
case $f = \exp X, X \in \Omega^{0,2,r+1} (TM)$. The main steps
in the proof are done already in \cite{ES}. We recall them.
Let $I=[0,1]$ and $ i_t:M \rightarrow I \times M$ the embedding 
$i_t(x)=(t,x)$   . We put on $I \times M$  the product metric
$ \left( \begin{array}{cc}
  1 & 0 \\     
  0 & g                  
  \end{array} \right)$    .

\me
\no
{\bf Lemma 3.9.} {\it 
For every $q \ge 0$ there exists a linear bounded mapping
\[ K: {}^{b,m} \Omega^{q+1} (I \times M) \rightarrow 
   {}^{b,m}\Omega^q(M)\]
such that $dK + Kd = i^*_1 - i^*_0$.
}

\me
\no
This is Lemma 3.1 of \cite{ES}.  \hfill $\Box$

\me 
\no   
{\bf Lemma 3.10.} {\it 
Let $f,h : M \rightarrow N$ be $C^1$--mappings and
$F: I \times M \rightarrow N$ a $C^1$--homotopy between $f$ and $h$.
Let 
\[ f^*, h^*: {}^{b,1}\Omega^q(N) \rightarrow {}^{b,1}\Omega^q(M), \quad
   F^*: {}^{b,1}\Omega^q(N) \rightarrow {}^{b,1}\Omega^q(I \times M) \]
be bounded. Then for $\Phi \in {}^{b,1}\Omega^q(N)$ 
\[ (h^*-f^*) \Phi = (dK + Kd) F^* \Phi .\]
}

\me
\no
This is lemma 3.2 of \cite{ES}. \hfill $\Box$

Hence we have to estimate $(dK+Kd) F^* \theta$ in our case $h = id$,
$f=\exp X$, $F=\exp tX$. This is done in theorem 3.2 of \cite{ES} and its 
proof, \cite{ES} p.154-158. The proof is rather involved. We conclude
\[ (\exp X)^* \theta - \theta \in \Omega^{1,2,r} . \]
Assume now $f=\exp X_n \circ \dots \circ \exp X_1$. A simple induction
now yields $f^* \theta - \theta \in \Omega^{1,2,r}$ (cf. \cite{ES} p. 160).
This finishes the proof of lemma 3.8 and hence of (3.1).
\hfill $\Box$

Now we consider (3.2) which is the differential of (3.1). From this it is 
clear that the expression (3.2) is in $\Omega^{2,2,r-1}$. But can we
prove more.
\be
  d \alpha \wedge f^* \theta + \alpha f^* d \theta - d \theta =
  d(\alpha - 1) \wedge f^* \theta + (\alpha-1) f^* d \theta +
  f^* d \theta - d \theta
\ee
and
\be
  d (\alpha-1) \wedge f^* \theta = d (\alpha-1) \wedge (f^* \theta - \theta)
  + d (\alpha-1) \wedge \theta ,
\ee
\be
  (\alpha-1) f^* d \theta = (\alpha-1)(f^* d \theta - d \theta) +
  (\alpha-1) d \theta .
\ee

Now we use that $(\alpha-1) \in \Omega^{0,2,r}$ and 
$f^* d \theta \in \Omega^{2,2,r}$ (according to the first part of the
proof) and $d \theta \in {}^{b,r} \Omega^2$. Application of the module
structure theorem 2.4 yields the assertion. This finishes the proof of
lemma 3.7. 
\hfill $\Box$

Define
\[ \Psi : {\cal F}^{r+1}_0 \timi {\cal D}^{r+1}_0 \rightarrow 
   \Omega^{1,2,r} \oplus \Omega^{2,2,r} , \]
\be
  \Psi (\alpha,f) := (\alpha f^* \theta - \theta, d \alpha \wedge
  f^* \theta + \alpha f^* \theta - d \theta).
\ee

\me 
\no             
{\bf Lemma 3.11.} {\it  The map   
$\Psi$ is  of class $k-r$. }

\me
\no
We omit the considerations and estimates, refer to \cite{EF2} and discuss in the
sequel only the special case of $D \Psi |_{(1, id)} $.
\hfill $\Box$

\me 
\no             
{\bf Lemma 3.12.} {\it     
Let 
$A: \Omega^{0,2,r+1} \oplus \Omega^{0,2,r+1}(TM) \rightarrow 
\Omega^{1,2,r} \oplus \Omega^{2,2,r} $
be defined by
\[ A(u, X) := ( u \cdot \theta + L_X \theta , d(u \cdot \theta) + 
   d (i_X d \theta )) . \]                        
Then                                             
\[ D  \Psi |_{(1, id)} (u, X) = A(u, X) . \]
}

\me
\no
{\bf Proof.}  From the facts that
${}^{b,r+1}|\theta| < \infty , u \in \Omega^{0,2,r+1}$ and
$X \in \Omega^{0,2,r+1} (TM)$  follows immediately that 
$u \cdot \theta , L_X \theta \in \Omega^{1,2,r}$ and
$d(u \cdot \theta), d(i_X d \theta) \in \Omega^{2,2,r}$.
Considering 
$\frac{d}{dt} \Psi (1+t \cdot u, id)|_{t=0},$ and 
$\frac{d}{d\tau} \Psi (1, \exp \tau X)|_{\tau=0}$
yields the desired result.
\hfill $\Box$

Define 
\[ B: \Omega^{1,2,r} \oplus \Omega^{2,2,r} \longrightarrow
   \Omega^{2,2,r-1} \oplus \Omega^{3,2,r-1} \]
by 
\[ B(\rho, \sigma) := (d \rho - \sigma, d \sigma) . \]
Clearly $BA = 0$.

\me
\no
{\bf Lemma 3.13.} {\it Let the adjoints $A^*,B^*$ be defined with respect
to the $L_2$ scalar product of forms. Then for
$\Box := AA^* + BB^*$, we have
\[ \Box (\rho, \sigma) =  
   (\bigtriangleup \rho + \rho , \bigtriangleup \sigma + \sigma ) . \]
}   

\me              
\no
We refer to lemma 8.3.2 of \cite{Om}. 
\hfill $\Box$
 
\me
\no
{\bf Corollary 3.14.} {\it
The operator $\Box$ is extendable to any Sobolev space of order $\le  k$
and
\[ \Box :  \Omega^{1,2,r} \oplus \Omega^{2,2,r} \longrightarrow
   \Omega^{2,2,r-2} \oplus \Omega^{3,2,r-2} \]
is injective, surjective and bounded.}
\hfill $\Box$

\me
\no
Consider $ker B \subset \Omega^{1,2,r} \oplus \Omega^{2,2,r}$. Then $\Psi$
maps ${\cal F}^{r+1} \timi {\cal D}^{r+1}_0 $ into
$\Omega^{1,2,r} \oplus \Omega^{2,2,r}$. The following is immediately
clear from the definitions. 

\me
\no
{\bf Lemma 3.15.} {\it
$im \Psi \subseteq ker B $ .}
\hfill $\Box$

\me 
\no
{\bf Proposition 3.16.} {\it
Assume that $\inf \sigma_e (\bigtriangleup_1 
|_{(ker \bigtriangleup_1)^\perp})>0$. Then the operator 
\[ D \Psi |_{(1,id)} = A: \Omega^{0,2,r+1} \oplus \Omega^{0,2,r+1} (TM)
   \longrightarrow T_{(0,0)} ker B = ker B \]
is surjective.}

\me
\no
{\bf Proof.} Consider
\[ \Omega^{0,2,r+1} \oplus \Omega^{2,2,r+1} (TM) 
   \begin{array}{c}                        
      A \\[-2ex] \longrightarrow \\[-2ex] {}
   \end{array}  
   \Omega^{1,2,r} \oplus \Omega^{2,2,r} 
   \begin{array}{c}
      B \\[-2ex] \longrightarrow \\[-2ex] {}
   \end{array}  
   \Omega^{2,2,r-1} \oplus \Omega^{3,2,r-1} . \]
This is an elliptic complex. Hence
\[ \Omega^{1,2,r} \oplus \Omega^{2,2,r} = ker \Box \oplus
   \overline{im A} \oplus \overline{im B^*} 
   = \overline{im A} \oplus \overline{in B^*} , \]
where the summands are $L_2$--orthogonal and the completion is taken with
respect to $|\en|_{2,r}$. Moreover, it is easy to see that 
$ker B \subseteq \overline{im A}$. Hence we are done if we can show that
\[ \overline{im A} = im A . \]
Now it is a well known fact from elementary functional analysis that $A$
is closed if and only if $im A A^*$ is closed. A longer calculation
yields
\be
  A A^* (\rho, \sigma) = (d \delta \rho + \rho + \delta \sigma, 
  d \delta \sigma + d \rho ) .
\ee
Hence $im A A^*$ is closed if and only if the operators
\[ (\rho, \sigma) \longrightarrow d \delta \rho + \rho + \delta \sigma \]
and
\[ (\rho, \sigma) \longrightarrow d \delta \sigma + d \rho \]
have closed image, respectively. Now a careful analysis shows that this
is the case if $im \bigtriangleup_1$ is closed. The latter is
equivalent to $\inf \sigma_e (\bigtriangleup_1 
|_{(ker \bigtriangleup_1)^\perp})>0$. We refer to \cite{EF2} for details. This
finishes the proof of proposition 3.16.
\hfill $\Box$

A series of shifting arguments yield the same result at any other point 
$(\alpha,f)$, i.e. 
$D \Psi |_{(\alpha,f)} $ is surjective. Hence $\Psi$ is a submersion and 
\[ {\cal D}^{r+1}_{\theta,0} = \Psi^{-1} (0,0) \]
is closed submanifold. The proof of proposition 3.6 is finished.
\hfill $\Box$

\me
\no
{\bf Corollary 3.17.} {\it
Assume $(M^{2n+1},g,\theta)$ satisfying $(I), (B_\infty), 
\sup\limits_{x \in M} |\nabla^i \theta|_x < \infty$ for all $i$ and 
$\inf \sigma_e (\bigtriangleup_1|_{(ker \bigtriangleup_1)^\perp}) > 0$.
Set ${\cal D}^\infty_{\theta,0} := 
\lim\limits_{\begin{array}{c} \leftarrow \\[-2ex] r \end{array}}
{\cal D}^{r+1}_{\theta,0}$.
Then
\[ \left\{ {\cal D}^{\infty}_{\theta,0}, {\cal D}^{r+1}_{\theta,0} 
   | r+1 \ge \frac{2n+1}{2}+2 \right\} \]
is an ILH Lie group. }
\hfill $\Box$

\setcounter{equation}{0}
 
\section{Contact Transformations of the Restricted Cotangent Bundle
        $T^* M \setminus 0$}

The most important example for us of contact manifolds of bounded 
geometry is the cotangent sphere bundle $S=S(T^*M) \cong 
((T^*M) \setminus 0) / {\bf R}_+$. We consider the Sasaki metric on $T^*M$.
Let $\pi : T^*M \rightarrow M$ be the projection and $K$ the connection
map of the Levi--Civita connection in the cotangent bundle. Then the
Sasaki metric is defined by
\[ g_{T^*M}(X,Y) = g_M(\pi_*X,\pi_*Y) + g_M(KX,KY), \en X,Y \in TT^*M . \]
Let $g_S := g_{T^*M}|_{S(T^*M)}$ be the restriction of the Sasaki metric to
the cosphere bundle.

\me
\no
{\bf Lemma 4.1.} {\it 
If $(M,g_M)$ satisfies $(I), (B_{k+1}), 0 \le k \le \infty$
fixed, then $(S(T^*M),g_S)$ satisfies $(I), (B_k)$.
 }

\me
\no
We refer to \cite{ES}, p. 165 for the proof.
\hfill $\Box$

\me
\no         
Let $\theta$ be the canonical one form on $T^*M$, i. e. for
$X \in T_zT^*M$, $\theta(X) := z(\pi_*X)$. Then
$\theta_s = i^* \theta$ is a contact form on 
$S(T^*M),$ where $ i: S(T^*M) \rightarrow T^*M$ is the inclusion. As pointed
out in \cite{Bl1}, $g_S$ and $\theta_S$ are not associated but this is true for
$g'_S := \frac{1}{4} g_S$ and $\theta'_S := \frac{1}{2} \theta_S$.
In \cite{Bl1}, p. 132--135 the Reeb vector field $\xi$, the (1,1) tensor field
$\varphi$ and the covariant derivatives are explicitly calculated.
Denote for the sake of simplicity the new $\theta'_S$ from now on by
$\theta = \theta'_S$ and $g'_S$ by $g_S$ .
           
\me        
\no        
{\bf Lemma 4.2.} {\it
Suppose $(M,g)$ with $(I), (B_{k+2})$. Then 
$\theta = \theta_s \in {}^{b,k+1} \Omega^1$, i. e.
$\sup\limits_{z \in S} |\nabla^i \theta|_{g_S,z} < \infty$, 
$0 \le i \le k+1$.
}          
           
\me        
\no        
{\bf Proof.} Start with $i=0$. Let $e_1, \dots, e_{2n-1}$ be an
orthonormal basis in $T_{z_0} S$ such that $e_{2n-1}=\xi$.
Then       
\[ |\Omega|^2_{g_S,z_0} = \sum\limits^{2n-1}_{i=1} 
   \Omega_{z_0}(e_i)^2 = 1 . \]
It is well known that on a contact Riemannian manifold the integral curves
of the Reeb vector field $\xi$ are geodesics (cf. \cite{Bl1} p. 54), i. e.
$\nabla_\xi \xi=0$. Fix at $z_0$ the orthonormal basis 
$e_1, \dots, e_{2n-1}, e_{2n-1}=\xi$. According to 
$\theta (X) = g(\xi, X)$, $\theta$ is the covariant form of $\xi$.
Hence $|\nabla^\nu \theta|=|\nabla^\nu \xi|$. 
\[ |\nabla \theta|^2 = |\nabla \xi|^2 = \sum\limits^{2n-1}_{i=1} 
   |\nabla_{e_i} \cdot \xi|^2 = \sum\limits^{2n-2}_{i=1} 
   |\nabla_{e_i} \cdot \xi|^2 . \]
According to \cite{Bl1}, p. 133--135, formulas (3)--(8), 
\bea      
  |\nabla_{e_i} \cdot \xi| & \le & C'_1 \cdot |\varphi| + C'_2 |R^{g_S}|~ ,~~ 
  |\nabla \xi| \le C_1 |\varphi| + C_2 |R^{g_S}| , \\
  |\nabla^\nu \theta| & \le & C_{1,\nu} |\nabla^{\nu-1} \varphi| +
C_{2,\nu} 
  |\nabla^{\nu-1} R^{g_S}| , \\
|\nabla^\mu \varphi| & \le & D_{1,\mu} |\nabla^{\mu-1} \varphi| +
D_{2,\mu} 
  |\nabla^{\mu-1} R^{g_S}|
\eea
which yields together with 4.1 the assertion.
\hfill $\Box$

\me
\no
{\bf Theorem 4.3.} {\it
Suppose $(M^n,g)$ with $(I)$ and $(B_{k+1})$, 
$k \ge r+1 > \frac{2n-1}{2}+2$, and 
$\inf \sigma_e                                             
(\bigtriangleup_1(g_S))|_{(ker \bigtriangleup_1(g_S))^\perp})>0$.  
Then ${\cal D}^{r+1}_{\theta,0} (S(T^*M),g_S)$ is well defined and a
$C^{k-r}$ submanifold of ${\cal F}^{r+1}_0 \timi {\cal D}^{r+1}_0$.
}
           
\me         
\no
This follows immediately from theorem 3.6. 
\hfill $\Box$
 
\me
\no 
{\bf Corollary 4.4.} {\it
Suppose $(M^n,g)$ with $(I)$ and $(B_\infty)$ and 
$\inf \sigma_e                                             
(\bigtriangleup_1(g_S))|_{(ker \bigtriangleup_1(g_S))^\perp})>0$.  
Set 
${\cal D}^\infty_{\theta,0} := \lim\limits_{\leftarrow}  
{\cal D}^{r+1}_{\theta,0}$.
Then
$\left\{ {\cal D}^\infty_{\theta,0} , {\cal D}^{r+1}_{\theta,0} |
r+1 > \frac{2n-1}{2}+2 \right\}$ is an ILH Lie group.
}
\hfill $\Box$

\me
\no
For our later applications we must rewrite 4.3 and 4.4 by rewriting
${\cal D}^{r+1}_{\theta, 0} (S(T^*M),g_S)$
 in an isomorphic version as ( writing ${\cal D}^{r+1}_{\theta,0} (S) $ for
short) 
\[ {\cal D}^{r+1}_{\theta,0} (S) =
   \left\{ (f,\beta) \in {\cal D}^{r+1}_0 (S) \timi
   {\cal F}^{r+1}_0(S) | f^* \theta = \beta \theta \right\} , \]
where $(\alpha,f) \mapsto (f,\alpha^{-1})$
is the canonical isomorphism (w.r.t. $\alpha$ antiisomorphism). The
''Lie algebra'' of ${\cal D}^{r+1}_{\theta,0} (S(T^*M))$ is 
\[ 
 {\bf d} ^{r+1}_{\theta,0}(S) = \left\{ (X,u) \in \Omega^{0,2,r+1} (TS)
   \timi \Omega^{0,2,r+1} (S) | L_X \theta = u \cdot \theta \right\}\]
with
\be
  [(X,u),(Y,v)]=([X,Y],X(v)-Y(u)) .
\ee
From the last equation it is clear that it isn't a Lie algebra since the
bracket decreases the Sobolev index. It is only the tangent space at
$(id,1)$.

Define now a map $\Phi$ from ${\cal D}^{r+1}_{\theta,0}(S)$ into the
homogeneous of degree one $C^2$ diffeomorphisms $\tilde{f}$ of
$T^* \setminus 0$ satisfying $\tilde{f}^* \theta = \theta$ (cf. 4.7
below). Given $(f,\beta) \in {\cal D}^{r+1}_{\theta,0}(S)$, we define
$\tilde{f}=\Phi(f,\beta)$ by
\be
  \tilde{f} (z) := \frac{f(\frac{z}{|z|})\cdot|z|}{\beta(\frac{z}{|z|})} ,~~~
  z \in T^*M \setminus 0 .
\ee
From $\tilde{f} \in im \Phi$ we can reproduce 
$(f,\beta)=\Phi^{-1}(\tilde{f})$ by
\be
  f(\frac{z}{|z|}) = \frac{\tilde{f}(z)}{|\tilde{f}(z)|} =
  \frac{\tilde{f}(\frac{z}{|z|})}{|\tilde{f}(\frac{z}{|z|})|} =
  \frac{\tilde{f}|_S(\frac{z}{|z|})}{|\tilde{f}|_S(\frac{z}{|z|})|}  
  \en , \quad 
  \beta(\frac{z}{|z|}) = \frac{|z|}{|\tilde{f}(z)|} =
  \frac{|\frac{z}{|z|}|}{|\tilde{f}(\frac{z}{|z|})|} =
  \frac{1}{|\tilde{f}|_S(\frac{z}{|z|})|}  .
\ee                                      
For our applications below we recall two lemmas from \cite{RS}.

\me
\no      
{\bf Lemma 4.5.} {\it   
Let $f:T^*M \setminus 0 \rightarrow T^*M \setminus 0$ be a diffeomorphism.
Then the following conditions are equivalent

\no
{\bf a)} $f^* \theta = \theta$

\no
{\bf b)} $f$ is symplectic (i. e. $f^*\omega=\omega, \omega=-d\theta$)
and homogeneous of degree one.
}                                  
\hfill $\Box$

\me
\no      
{\bf Lemma 4.6.} {\it   
{\bf a)}
Let $H: T^*M \setminus 0 \rightarrow R$ be homogeneous of degree $\nu$.
Then the Hamiltonian vector field $X_H$ is homogeneous of degree
$\nu-1$ and $\theta(X_H)=H$.

\no
{\bf b)}
A vector field $X$ on $T^*M \setminus 0$ is homogeneous of degree zero
if and only if its flow is homogeneous of degree one.

\no
{\bf c)}
$L_X \theta = 0$ if and only if $X$ is globally Hamiltonian, homogeneous of
degree zero with Hamiltonian $\theta(X)$ homogeneous of degree one.
}

\me
\no                                            
{\bf Proof.}                                                      
For later use we recall the proof of c). Let $L_X \theta = 0$, $F_t$ the 
flow of $X$. Then $F_t^* \theta = \theta$. Lemma 4.5 implies that $F_t$
is symplectic and homogeneous of degree one and hence,
according to b), $X$  is of degree zero.
$0=L_X\theta=i_Xd\theta+di_X\theta$ yields                        
$i_X\omega=d\theta(X)$, i. e. $X=X_H$ with $H=\theta(X)$ homogeneous               
of degree one. The converse implication follows from a).
\hfill $\Box$                                                                   

\me                                                                           
\no
Define now                                                                  
\[ {\cal D}^{r+1}_{\theta,0} (T^*M \setminus 0) := \left\{ 
   \tilde{f} : T^*M \setminus 0                    
   \begin{array}{c}                                                         
       \cong \\[-2ex] \longrightarrow                                        
   \end{array}                                                              
   T^*M \setminus 0 | \tilde{f} = \Phi(f,\beta), \en 
   (f,\beta) \in {\cal D}^{r+1}_{\theta,0} (S) 
   \right\} \]                                    
Recall our assumptions, $(M^n,g)$ with $(I),(B_{k+2}), 
k \ge r+1 > \frac{2n-1}{2}+2$. This implies $\omega \in 
{}^{b,r+1} \Omega^2(T^*M)$. Additionally we have in the case of $T^*M$
 that $\omega$ is strongly nondegenerate, i. e.            
$\inf\limits_{z \in T^*M}|\omega|^2_z>0$.         
                                                  
It follows immediately from the definition (4.5) that $\tilde{f}$
is a $C^2$ diffeomorphism. Thus we get a 1--1 mapping between
${\cal D}^{r+1}_{\theta,0}(S)$ and                 
${\cal D}^{r+1}_{\theta,0}(T^*M \setminus 0)$. We endow 
${\cal D}^{r+1}_{\theta,0}(T^*M \setminus 0)$   with the topology and
differential structure of
${\cal D}^{r+1}_{\theta,0}(S)$ such that  $\Phi$ becomes a diffeomorphism. 
Evidently, $\Phi(id,1)=id_{T^*M \setminus 0}$. Our next aim is to
describe properties of $\Phi$, of the elements of  
${\cal D}^{r+1}_{\theta,0}(T^*M \setminus 0)$ and of 
$T_{id}{\cal D}^{r+1}_{\theta,0}(T^*M \setminus 0)$.
                                                   
\me                                                
\no                                                
{\bf Proposition 4.7.} {\it                        
                                                   
\no                                                
{\bf a)}                                           
$\Phi$ is an isomorphism of groups.                
                                                   
\no                                                
{\bf b)}                                           
Each $\tilde{f} \in {\cal D}^{r+1}_{\theta,0}(T^*M \setminus 0)$ satisfies
$\tilde{f}^* \theta = \theta$.                     
                                                   
\no                                                
{\bf c)}                                           
Let $(X,h) \in T_{(id,1)} {\cal D}^{r+1}_{\theta,0} (S)$. Then
$\Phi_{*(id,1)} (X,u) = X_H$  where $H$ is the Hamiltonian 
$H(z)= \theta_S (X_{\frac{z}{|z|}}) \cdot |z|$, i. e. 
\be
  H = \theta_S (X \circ \pi) / f_S, \en            
  f_S (z) = \frac{1}{|z|} .
\ee                                                
In particular, $H$ is homogeneous of degree one, $X_H$ is
homogeneous of degree 0 and $H= \theta(X_H)$.      

\no                                                
{\bf d)}
If $Y=Y_H \in T_{id}{\cal D}^{r+1}_{\theta,0}(T^*M \setminus 0)$ 
then $Y$ projects to $X=\pi_*Y$ tangentially to $S$ and 
$\Phi^{-1}_* (Y) \equiv (\Phi^{-1})_*(Y) = (\pi_* Y,u)$, where
$u(\frac{z}{|z|})= \left\{ \frac{1}{|z|} ~, H(z) \right\} 
\cdot |z|$, i. e.                                   
\be
   u \circ \pi = \left\{ f_S, H \right\} / f_S       
\ee
Here \{ , \} the Poisson bracket on $T^*M$.           

\no                                                    
{\bf e)}
$\Phi_*$ preserves the Lie brackets.                     

\no                            
{\bf f)}                        
$T_{id}{\cal D}^{r+1}_{\theta,0}(T^*M \setminus 0)$ coincides with the space
\bea                                                    
   {\cal H}^0_\theta \Omega^{0,2,r+1} (T(T^*M \setminus 0)) & = &
   \Big\{ Y \en | \en Y ~is~ a ~ C^2 \en \mbox{vector field on} \en T^*M
\setminus  0, \en
   L_Y \theta = 0 \nonumber \\           
   &{}& \mbox{and} \en Y|_S \in                          
   \Omega^{0,2,r+1} (i^*(T(T^*M \setminus 0))), \en 
   i: S \rightarrow T^*M \setminus 0                    
   \Big\} . \nonumber                                    
\eea
                       
\no                 
{\bf g)}           
${\cal H}^0_\theta \Omega^{0,2,r+1} (T(T^*M \setminus 0))$ is isomorphic
to the space                
\bea                
  {\cal H}^1 \Omega^{0,2,r+2} (T^*M \setminus 0) & = &
  \Big\{ h ~| h ~is ~a~ C^3 ~function ~on ~
  T^*M \setminus 0, \nonumber \\
  &{}& h \en \mbox{is homogeneous of degree 1 and} \nonumber \\
  &{}& i^*h=h|_S \in \Omega^{0,2,r+2} (S) \Big\}. \nonumber
\eea

\me
\no
{\bf h)}
${\cal H}^1 \Omega^{0,2,r+2} (T^*M \setminus 0)$ is isomorphic to 
$\Omega^{0,2,r+2} (S)$.
}

\me
\no
{\bf Proof.}
For a) we refer to \cite{RS}, p. 97. Let $\tilde{f} \in 
{\cal D}^{r+1}_{\theta,0}(T^*M \setminus 0)$. Recall
$\theta_S = \theta|_S =i^* \theta$, $\theta$ on
$T^*M \setminus 0$ the canonical one form, $\pi: T^*M \setminus 0 \rightarrow
S(T^*M)=S$ the projection, $\pi(z)=\frac{z}{|z|}$. Then 
$\tilde{f}^* \theta = \left( \frac{i \circ f \circ \pi}
{(\beta \circ \pi) \cdot f_S} \right)^* \theta = 
\frac{\pi^*f^*\theta_S}{(\beta \circ \pi) \cdot f_S} =
\frac{(\beta \circ \pi)\pi^*\theta_S}{(\beta \circ \pi) f_S} =
\theta$. This proves b). We conclude, according to 4.5.b) that $\tilde{f}$
is symplectic and homogeneous of degree one. Claims c), d) and e) are again
simple calculations, performed in \cite{RS} p. 97,98. For claim f) we use
the fact that 
$X_H$ is homogeneous of degree 0 to assure that $Y_H$ projects to
$\pi_* Y_H$ tangentially to $S$. Let
$Y=X_H=\Phi_*(X,u) \in T_{id} {\cal D}^{r+1}_{\theta,0} (T^*M \setminus 0)
\equiv \Phi_* T_{id}{\cal D}^{r+1}_{\theta,0}(S)$. Assuming for a moment
$Y=X_H \in C^2$, we conclude from 4.6.c) that $L_X \theta = 0$. Hence we
have only to show that $Y=X_H$ is $C^2$ and $Y|_S \in 
\Omega^{0,2,r+1} (i^*(T(T^*M \setminus 0)))$. The latter would imply that
$Y|_S \in C^2$ (according to the Sobolev embedding theorem), hence
$Y \in C^2$ since $Y$ is homogeneous of degree zero. Assuming 
$Y|_S = X_H |_S \in \Omega^{0,2,r+1} (i^*(T(T^*M \setminus 0)))$, we
have $T_{id} {\cal D}^{r+1}_0 (T^*M \setminus 0) \subseteq 
{\cal H}^0_\theta \Omega^{0,2,r+1}(T(T^*M \setminus 0))$. Consider
$\supseteq$. Let $Y \in {\cal H}^0_\theta \Omega^{0,2,r+1}(T(T^*M 
\setminus 0))$. Then, according to 4.6.c) $Y=Y_H$ for some $H$ and $Y$ has 
degree zero homogeneity. Hence it projects to $T S$, 
$\pi_* Y  = X$. Define $u$ by (4.8). Assuming for a moment 
$X \in \Omega^{0,2,r+1}(T S)$ and $u \in \Omega^{0,2,r+1}(S)$, we 
see by an easy calculation $\Phi_*(X,u)=Y_H$, i. e. $\supseteq$ would be
proved. Hence there remains to show 
1. $Y=X_H \in T_{id} {\cal D}^{r+1}_{\theta,0}(T^*M \setminus 0)$
implies $Y|_S \in \Omega^{0,2,r+1}(i^*(T(T^*M \setminus 0)))$,
2. $Y \in {\cal H}^0_\theta \Omega^{0,2,r+1}$, $Y=Y_H$ then
$X= \pi_* Y \in \Omega^{0,2,r+1}(TS)$ and 
$u \in \Omega^{0,2,r+1}(S)$. 

Lets begin with the first assertion $Y=X_H=\Phi_*(X,u)$.
We introduce local coordinates in $T(T^*M \setminus 0)$, say
$(x,\xi,\psi_1,\psi_2)$, $x$ coordinates im $M$, $\xi$ in 
$T^*M \setminus 0$, $\psi_1, \psi_2$ over them in $T(T^*M \setminus 0)$
with projections $\pi_1, \pi_2$. Then according to \cite{RS}, p. 99 
we have
\be
  X_H|_S(x,\xi) = \Phi_*(X,u)|_S(x,\xi) = 
  \Big( x,\xi,(X)_1(x,\xi),(X)_2(x,\xi)-u(x,\xi) \cdot \xi \Big) ,
\ee
where $( ~~ )_i$ are the components of $X$ belonging to $im \pi_i$.
Using a uniformly locally finite cover of $S(T^* M \setminus 0)$,
(4.9) and $X \in \Omega^{0,2,r+1}(TS), u \in \Omega^{0,2,r+1}(S)$,
we conclude that $Y|_S=X_H|_S \in \Omega^{0,2,r+1}
(i^*(T(T^*M \setminus 0)))$.   
                         
To prove 2., we assume $Y=Y_H \in {\cal H}^0_\theta \Omega^{0,2,r+1}(T(T^*M
\setminus 0))$. By assumption $Y|_S \in \Omega^{0,2,r+1}
(i^*(T(T^*M \setminus 0)))$. Set $X=\pi_*Y=\pi_*Y|_S$. Denote by
$\nabla^S$ the Levi-Civita connection of $(S(T^*M),g_S)$.
Then by choice of local orthonormal bases $e_1, \dots, e_{2n-1},e_{2n},
e_{2n} \perp S $, we see immediately for $i \le r+1$
\be
  |X|_{g_S} \le |Y|_S| , \quad |(\nabla^S)^i X | \le 
  |\nabla^i (Y|_S )| ,
\ee
which implies $X \in \Omega^{0,2,r+1}(TS)$. Write as in (4.9)
\be
  X(x,\xi)=(x,\xi,(X)_1(x,\xi), (X)_2(x,\xi)) .
\ee
Then, locally,
\be
  (X_H|_S-X)(x,\xi)=(0,0,0,-u(x,\xi) \cdot \xi) ,
\ee
which immediately implies $u \in \Omega^{0,2,r+1}(S).$

\no
This finishes the proof of f).

Consider g) and the map $Y=X_H \mapsto H=\theta(X_H)$ given by
4.6.c. We must prove that $H \in C^3$ and $H|_S=i^*H \in 
\Omega^{0,2,r+2}(S)$. The latter will already imply $H \in C^3$.
We immediately obtain from (4.7) that 
$H|_S \in \Omega^{0,2,r+1}(S)$ since $H=\theta(X_H)$,
$H|_S=\theta(X_H|_S)$, $\theta \in {}^{b,r+1}\Omega^1$ and
$X_H|_S \in \Omega^{0,2,r+1}(TS)$. Hence $H \in C^2$. The main
point is that $H|_S$ has even Sobolev order $r+2$. Denote again by
$d_S, \nabla^S$ the operators $d, \nabla$ on $S$. We have
\be
  |d_S(H|_S)| \le |(dH)|_S|, |(\nabla^S)^i(H|_S)|
  \le |(\nabla^iH)|_S| .
\ee
For nonsmooth objects we have (as usual) to understand this in the 
distributional sense. We always have to do with regular distributions.
That $H|_S\in \Omega^{0,2,r+2}(S)$ would be proved if we could show that
$d_S H|_S \in \Omega^{1,2,r+1}(S)$. According to (4.13) this
would be done if we could show
$(dH)|_S, ,(\nabla dH)|_S \dots (\nabla^{r+1} dH)|_S  $ are
square integrable on $S$. But $(dH)|_S = (i_{X_H} \omega)|_S
=i_{X_H|_{S}} \omega$. Furthermore $\omega \in {}^{b,r+2}\Omega^2$, $\omega$ 
is strongly nondegenerate and $X_H|_{S} \in \Omega^{0,2,r+1}
(i^*(T(T^*M \setminus 0)))$ just imply that $(dH)|_S , \dots , 
(\nabla^{r+1}dH)|_S$ are square integrable on $M$. This is Lemma 3.6, 3.7 in
\cite{ES}. We obtained that $H \in {\cal H}^1 \Omega^{0,2,r+2}(T^*M \setminus
0)$, and from 4.6 it follows that the map $X_H \mapsto \theta(X_H)
\mapsto X_{\theta(X_H)}$ equals to the identity.
          
Now let $h \in {\cal H}^1 \Omega^{0,2,r+2}(T^*M \setminus 0)$. The function 
$h$ defines a global Hamiltonian vector field $X_h$, homogeneous of degree
zero, satisfying $L_{X_h} \theta =0$. Moreover $X_h$ is $C^2$ and solves the
equation 
$i_{X_h} \omega = dh$. We have to assure that $X_h|_S \in \Omega^{0,2,r+1}
(i^*(T(T^*M \setminus 0)))$. From $h|_S \in \Omega^{0,2,r+2}(S)$
we conclude that $|h|_S, |d_S h|$ and $ |(\nabla^S)^i d_S h|  $ are square
integrable on $S$. But for $i \ge 1$
\[ |(\nabla^ih)|_S \le |(\nabla^S)^{i-1} (h|_S)|
   + |(\nabla^S)^i (h|_S)| , \]
since $h$ is homogeneous of degree one. We obtain that
$|(dh)|_S, \dots , |(\nabla^{r+1}dh)|_S$ are square
integrable on $S$, i. e. the right hand side of 
$i_{X_h} \omega_S = dh|_S$ is an element of 
$\Omega^{1,2,r+1}(S)$ (with values in the conormal bundle of
$S$). Then $X_h  \in \Omega^{0,2,r+1} (i^*(T(T^*M \setminus 0)))$ and 
altogether $X_h \in {\cal H}^0_\theta \Omega^{0,2,r+1}  
(T(T^*M \setminus 0))$. According to 4.6, the map
$h \mapsto X_h \mapsto \theta(X_h)$ coincides with $id$.
This finishes th proof of claim g).
            
Concerning claim h), the isomorphism is given by $h \in {\cal H}^1 
\Omega^{0,2,r+2} (T^*M \setminus 0) \mapsto h|_S$ .
This map is well defined, according to g). It is injective since
$h$ is homogeneous of degree one. It is surjective because for 
$u \in \Omega^{0,2,r+2}(S)$ let $h_u$ be its extension homogeneous of
degree one . Then $h_u \in {\cal H}^1 \Omega^{0,2,r+2} (T^*M
\setminus 0)$  and $h_u|_{S}=u$. 
\hfill $\Box$
         
\me         
\no      
{\bf Remark.} We constructed a topological isomorphism
\be       
  F: {\bf d}^{r+1}_{\theta,0}(S) \longrightarrow \Omega^{0,2,r+2} (S) .
\ee       
This isomorphism is topological since all constructed maps in Proposition 4.7
are  norm continuous. Here we essentially use Lemma 3.6, 3.7 of \cite{ES}. The
isomorphism (4.14) will be very important in constructing local charts on
${\cal D}^{r+1}_{\theta,0}(S)$. Proposition 4.7. justifies to denote
${\cal D}^{r+1}_{\theta,0}(S)$ and 
${\cal D}^{r+1}_{\theta,0}(T^*M \setminus 0)$ by the same symbol. We can
understand ${\cal D}^{r+1}_{\theta,0}(T^*M \setminus 0)$ as a Hilbert manifold
and a topological group with 
$T_{id} {\cal D}^{r+1}_{\theta,0}(T^*M \setminus 0) \cong 
\Omega^{0,2,r+2}(S)$.
\hfill $\Box$
            
\me        
Summarizing our results, we obtain in the case of $(B_\infty)$ the 
following : 
           
\me        
\no        
{\bf Theorem 4.8.} {\it
Suppose $(M^n,g)$ with $(I), (B_\infty)$ and 
$\inf \sigma_e(\bigtriangleup_1(g_S)|
_{ker \bigtriangleup_1(g_S)^\perp}) >0$.
Set ${\cal D}^\infty_{\theta,0}(T^*M \setminus 0) = 
\lim\limits_{\leftarrow} {\cal D}^{r+1}_{\theta,0}(T^*M \setminus 0)$.
Then         
\[ \Big\{ {\cal D}^\infty_{\theta,0}(T^*M \setminus 0)~,~ 
   {\cal D}^{r+1}_{\theta,0}(T^*M \setminus 0) ~\Big| ~
   r+1 > \frac{2n-2}{2}+2 \Big\} \]
is an ILH Lie group.
}             
\hfill $\Box$ 
                 
\setcounter{equation}{0}
               
\section{Pseudodifferential and Fourier Integral Operators on Open
Manifolds}

Pseudodifferential ($\Psi$DO) and Fourier integral operators (FIO) are
well defined for any manifold, open or closed. But on open manifolds the
spaces of these operators don't have any reasonable structure. Moreover, many
theorems for
$\Psi$DOs or FIOs on closed manifolds become wrong or don't make any sense in
the open case, e. g. certain mapping properties between Sobolev spaces of
functions are wrong. The situation rapidly changes if we restrict ourselves to
bounded geometry and adapt these operators to the bounded geometry. This
means, roughly speaking, that the family of local symbols together with their
derivatives should be uniformly bounded. For FIOs we additionally restrict
ourselves 
to comparatively smooth Lagrangian submanifolds $\Lambda$ of $T^*M \setminus
0 \times T^*M
\setminus 0$  and phase functions also adapted to the bounded geometry.

A good reference for $\Psi$DO's are \cite{Ko} and \cite{Sh}. Further results are in
preparation (cf. \cite{EK}). Since we restrict our applications to the
case where the Riemannian manifold $(M^n,g)$ satisfies the conditions of
bounded geometry $(I)$ and $(B_\infty)$, we assume these conditions from now
on. Moreover, we restrict ourselves to the scalar case, i. e. we consider
only  operators acting on functions.

We first recall two classical lemmas which play a key role in all
forthcoming constructions.

\me
\no    
{\bf Lemma 5.1.} {\it   
Assume $(M^n,g)$ with $(I)$ and $(B_\infty)$ ($(B_0)$ is   sufficient here), 
$\delta < \frac{r_{inj}}{2}$. Then there exists a uniformly locally
finite cover ${\cal U}=\{U_i\}_i$ of $M$ by geodesic $\delta$--balls.
}                                                
\hfill $\Box$

\me               
\no               
{\bf Lemma 5.2.} {\it   
Assume ${\cal U} = \{U_i\}_i$ as above (and $(B_\infty)$). Then there
exists an associated partition of unity $\{\psi_i\}_i$ such that 
\be             
  |\nabla^k \psi_i| \le C_k, \quad k=0,1,2, \dots, .
\ee             
}               
\hfill $\Box$  
               
\me            
We define now $U\Psi^{-\infty}(M)$ to be the set of all linear
operators $R:C^\infty_c(M) \rightarrow C^\infty_c(M)$ which have the
following properties.
              
\no
1. $R$ has Schwartz kernel ${\cal K}_R \in C^\infty(M \times M)$.
\hfill (5.2)
            
\no         
2. There exists a constant $C_R$ s. t. ${\cal K}_R (x,y) =0$
for $d(x,y) > C_R$. 
\hfill (5.3)
          
\no       
3. $\nabla^i_x \nabla^j_y {\cal K}_R$ is bounded for all $i$ and $j$.
\hfill (5.4)
        
\no    
It follows from the conditions $(I)$ and $(B_\infty)$ that for any point $m
\in M^n$ there exists a diffeomorphism $\Phi_m$ 
\[ M \supset B_\varepsilon(m) \begin{array}{c} \Phi_m \\[-2ex] \longrightarrow
   \\[-2ex] {} \end{array} B = B_\varepsilon(0) \subset {\bf R}^n \]
such that $\Phi_m$ induces bounded isomorphisms
\[ {}^{b,k}\Omega^0(B_\varepsilon(m)) \begin{array}{c} \cong \\[-2ex] 
\longrightarrow
   \\[-2ex] {} \end{array} {}^{b,k}\Omega^0(B) , \]  
\[ {}^{b,\infty}\Omega^0(B_\varepsilon(m)) \begin{array}{c} \cong \\[-2ex] 
\longrightarrow
   \\[-2ex] {} \end{array} {}^{b,\infty}\Omega^0(B)  \]  
with bounds independent of $m$. After fixing an orthonormal basis in
$T_mM$, $\Phi_m$ is essentially given by the exponential map.
              
\setcounter{equation}{4}              
We now define the class of uniform symbols for our pseudodifferential
operators as follows. Let $q \in {\bf R}$ and denote by ${\cal U} S^q (B)$ the
set of all families $\{a_m\}_{m \in M}$ with $a_m \in C^\infty (B \times {\bf
R}^n)$ and 
\be
  |\partial^\alpha_\xi \partial^\beta_x a_m (x,\xi)| \le 
  C_{\alpha,\beta} (1+|\xi|)^{q-|\alpha|} ,
\ee
where $C_{\alpha,\beta}$ is independent of $m$. Then $\{a_m\}_m$ defines
a family of operators
\[ a_m(x,D_x): C^\infty_c (B) \longrightarrow C^\infty (B)\]
by
\be
  a_m(x,D_x) u(x) := (2\pi)^{-n} \int\limits_{{\bf R}^n} \int\limits_B
  a_m (x,\xi) u(y) dy d\xi ,
\ee
where supp $u \subset B$.   

Define ${\cal U}\Psi^{-\infty}(B)$ as the set of all families
$\{R_m: C^\infty_c(B) \rightarrow C^\infty(B)\}_{m \in M}$
such that $R_m$ has Schwartz kernel ${\cal K} R_m \in C^\infty(B \times B)$
with 
\be
  |\partial^\alpha_x \partial^\beta_y {\cal K}_{R_m}(x,y)| \le
  C_{\alpha,\beta} ,
\ee
where $ C_{\alpha,\beta}$ is independent of $m$. Finally let
${\cal U}\Psi^q(B)$ be the set of all families
$\{A_m: C^\infty_c(B) \rightarrow C^\infty(B)\}_{m \in M}$
such that 
\[ A_m = a_m(x,D_x) + R_m , \]
$\{a_m\}_m \in {\cal U} S^q (B), 
\{R_m\}_m \in {\cal U} \Psi^{-\infty} (B)$.

Now we define the space ${\cal U} \Psi^q(M)$ of {\em uniform
pseudodifferential operators of order $q$}  on $M^n$ as follows:
A pseudodifferential operator $A$ on $M^n$ with Schwartz kernel ${\cal K}_A$
belongs to ${\cal U} \Psi^q(M)$ iff it satisfies the following conditions:

\no
1. There exists a constant $C_A>0$ s. t.
\be
  {\cal K}_A (x,y) = 0 \en \mbox{for} \en d(x,y)>C_A, x,y \in M .
\ee

\no
2. ${\cal K}_A$ is smooth outside the diagonal of $M \times M$.

\no
3. For any $\delta >0$ and $i,j$ there exists a constant
$C_{\delta,i,j}>0$ s. t. 
\be
  |\nabla^i_x \nabla^j_y {\cal K}_A (x,y)| \le C_{\delta,i,j}
  \en \mbox{for} \en d(x,y)>\delta .
\ee

\no 
4. If $A_m$ is defined by the following commutative diagram

\bi                                                   
                                                      
\[ C^\infty_C (B_\varepsilon(m))                      
   \quad                                              
   \begin{array}{c} A \\                              
   \Large\longrightarrow \\ {}                       
   \end{array}                                       
   \quad                                             
   C^\infty (B_\varepsilon(m)) \]                    
\vspace{-2ex}                                        
\[ \exp^*_m \Bigg\downarrow \hspace{4cm}             
   \Bigg\downarrow \exp^*_m \]                       
\vspace{-1ex}                                        
\[                                                   
   C^\infty_C(B)                                      
   \quad \quad                                        
   \begin{array}{c} A_m \\                            
   \Large\longrightarrow \\ {}                        
   \end{array}                                        
   \quad \quad                                        
   C^\infty(B) \]                                     
                                                      
\bi                                                   
\no                                                   
then the family $\{A_m\}_{m \in M}$ belongs to ${\cal U}\Psi^q(B)$.
                                                     
\me                                               
\no                                            
{\bf Remark.} {\it                          
We have  $~{\cal U}\Psi^{-\infty}(M) = \bigcap\limits_q 
{\cal U}\Psi^q(M)$. }                     
\hfill $\Box$                            

\me
\no
A convenient description for the elements $A \in {\cal U}\Psi^q(M)$
is given by

\me
\no
{\bf Proposition 5.3.} {\it
Assume $A \in {\cal U}\Psi^q(M),$ and $ \varepsilon>0$ arbitrary. Then there
exists a representation $A=A_1+A_2$

a. $A_1 \in {\cal U}\Psi^{-\infty}(M)$

b. $A_2 \in {\cal U}\Psi^q(M)$

c. ${\cal K}_{A_2(x,y)} =0$ for $d(x,y)>\varepsilon$

\no
i. e. up to smoothing operators in ${\cal U}\Psi^{-\infty}(M)$ we can
always assume that the support of ${\cal K}_A$ is arbitrary dense to 
the diagonal.
}

\me
\no
We refer to \cite{Ko}, p. 230/231 for the proof.
\hfill $\Box$

\me
For our applications we additionally restrict ourselves to {\em classical}
symbols and {\em classical}
$\Psi$DOs, i. e. we assume homogeneity in the $\xi$--variable on 
${\bf R}^n \setminus \{ 0 \}$ and an asymptotic expansion
\be
  a_m (x,\xi ) \sim \sum\limits^\infty_{j=0} a_{m,q-j}(x,\xi) 
\ee
such that $a_{m,q-j}(x,\xi )$ is positive homogeneous of degree $q-j$ in
$\xi$. Here $a_m (x, \xi ) \sim \sum\limits^\infty_{j=0} a_{m,q-j}(x,\xi)$
means
\be
  \left\{ (1-\chi(x,\xi)) \left( a_m(x,\xi) -
  \sum\limits^{k-1}_{j=1} a_{m,q-j}(x,\xi) \right) \right\}_m
  \in {\cal U} S^{q-k}(B)
\ee 
for all $k$ and $\chi (x, \xi)$ is compactly supported in the
$\xi$--direction with $\chi=1$ in a neighborhood of $B \times \{0\}$.

\me
\no
{\bf Remark.} {\it
Shubin \cite{Sh} calls such pseudodifferential operators $\Psi$DOs with
polyhomogeneous local symbols and writes ${\cal U}\Psi^q_{phg}(M)$ but we
omit the subscript phg and write simply ${\cal U}\Psi^q(M)$. 
}

\me
\no
We recall some mapping properties and refer to \cite{Ko} for the proofs. 

\me
\no
{\bf Proposition 5.4.} {\it 
Any $R \in {\cal U}\Psi^{-\infty}(M)$ defines continuous maps and
extensions as follows

$R: C^\infty_c (M) \longrightarrow C^\infty_c (M)$ , 
  
$R: C^\infty (M) \longrightarrow C^\infty (M)$ ,

$R: {\cal E}' \longrightarrow C^\infty_c (M)$ ,  
 
$R: {}^{b,\infty}\Omega^0 (M) \longrightarrow {}^{b,\infty}\Omega^0 (M)$  .
}
\hfill $\Box$

\me
\no
{\bf Proposition 5.5.} {\it
$A \in {\cal U}\Psi^q(M)$ defines linear continuous maps
\[ A : \Omega^{0,2,r} \longrightarrow \Omega^{0,2,r-q} \]
and                            
\[ A : {}^{b,\infty} \Omega^0 \longrightarrow  {}^{b,\infty} \Omega^0  .\]
}                             
\hfill $\Box$

\me
\no                              
Finally we have
                              
\me                               
\no
{\bf Proposition 5.6.} {\it
                                       
a. If $R_1,R_2 \in {\cal U} \Psi^{-\infty} (M)$ then 
$R_1 \circ R_2 \in {\cal U} \Psi^{-\infty} (M)$. 
                                             
b. If $R \in {\cal U} \Psi^{-\infty} (M)$ then 
$R^* \in {\cal U} \Psi^{-\infty} (M)$.  
                   
c. If $A \in {\cal U} \Psi^q (M)$ then 
$A^* \in {\cal U} \Psi^q (M)$  

d. If $A \in {\cal U} \Psi^{q_1} (M)$, $B \in {\cal U} \Psi^{q_2} (M)$ then 
$A \circ B \in {\cal U} \Psi^{q_1+q_2} (M)$ .
}

\me
\no
All proofs are performed locally. Using the uniform boundedness (5.5), one
gets these results for the formulas of the symbols of the adjoint
operators and the product (composition) of operators.
\hfill $\Box$

\me
\no
Finally we recall uniform ellipticity.  $A \in {\cal U} \Psi^q (M)$ 
is called {\em uniformly elliptic} if there exist constants $C_1,C_2,R > 0$,
independent of $m \in M$, such that
\[ C_1 |\xi|^q \le |a_m(x,\xi)| \le C_2 |\xi|^q \]
for all $|\xi|>R, x \in B, m \in M$. We denote this class of operators by
$E{\cal U}\Psi^q(M)$.

\me
\no
{\bf Remark.} {\it
Given any real number $s$, then there exists a uniformly elliptic operator in
$E{\cal U}\Psi^s(M)$, e. g. $(1+\bigtriangleup)^{\frac{s}{2}} \in 
E{\cal U}\Psi^s(M)$.         
}

\me
\no
{\bf Theorem 5.8.} {\it
Given $A \in E{\cal U}\Psi^q(M)$,  then there exists a parametrix, i. e. a
$P \in {\cal U}\Psi^{-q}(M)$ s. t.
\be
  P \circ A = I + R_1, \quad A \circ P = I + R_2, \quad
  R_1, R_2 \in  {\cal U}\Psi^{-\infty}(M) .
\ee
}

\no
The proof is performed locally by establishing explicit formulas for the
symbol  $\{p_m\}_m$. This is done as usual by calculation of the terms of the
asymptotic  expansion (5.10). Then one fits the local operators together by a
partition of unity. To assure $P \in {\cal U}\Psi^{-q}(M)$, one
essentially uses (5.1).
\hfill $\Box$

\me
\no
{\bf Remark.} {\it
In contrast to the case of compact manifolds, (5.12) does not mean
the invertibility of $A$ modulo compact operators. On open manifolds, 
the kernels ${\cal K}_{R_1}, {\cal K}_{R_2}$ are far from being square
integrable, i. e. $R_1, R_2$ are far from being compact operators. As 
a simple consequence, $P$ is far from being Fredholm (except in very
special cases).
}
\hfill $\Box$

\me
\no
Taking 5.3 into account, we define as in \cite{ARS1} a formal $\Psi$DO of order
$q$ as an element of ${\cal U}\Psi^q(M) / {\cal U}\Psi^{-\infty}(M)$.
Denote by $({\cal U}\Psi(M))_*$ the set of all invertible elements in 
${\cal U}\Psi(M)$,
\[ {\cal U}\Psi(M) := \bigcup\limits_q 
   {\cal U}\Psi^q(M) / {\cal U}\Psi^{-\infty}(M) = 
   \left( \bigcup\limits_q 
   {\cal U}\Psi^q(M) \right) / {\cal U}\Psi^{-\infty}(M) . \]
$({\cal U}\Psi(M))_*$ is a graded group under multiplication and non--empty
as 5.8 shows.

Quite similarly we define uniform Fourier integral operators 
${\cal U}F^q(M,C)$. A Fourier integral operator (FIO) on $M$ has 
essentially 3 ingredients

1. a family $a=\{a_m\}_m$ of local uniform symbols as above,

2. a conic Lagrangian submanifold (or homogeneous canonical relation)
$C \subset T^*M \setminus 0 \times T^*M \setminus 0$,

3. a family $\varphi=\{ \varphi_m\}_{m \in M}$ of phase functions.

We will make this precise. Recall that we now always assume that
$(M^n,g)$ satisfies the conditions $(I)$ and $(B_\infty)$. 
A homogeneous
canonical relation
$C$ is a closed submanifold (not in the sense of compactness)
$C \subset T^*M \setminus 0 \times T^*M \setminus 0$
which is conical, i. e. $((x,\xi),(y,\eta) \in C$ and $\tau>0$
imply $((x,\tau\xi),(y,\tau\eta) \in C$, and which is Lagrangian with
respect to the canonical symplectic form 
$\omega \ominus \omega = p^*_1 \omega - p^*_2 \omega$, where $p_j :T^*M
\times T^*M \longrightarrow T^*M , j=1,2$ are the projections and
$\omega = -d \theta$. 
A very important class of examples are the
graphs $\Gamma(\tilde{f})$ of contact transformations.
Let $f \in {\cal D}^{r+1}_{\theta,0} (S(T^*M))$ and  
$\tilde{f}=\Phi(f)$ the corresponding homogeneous diffeomorphism
$\tilde{f}:T^*M \setminus 0 \longrightarrow T^*M \setminus 0 $
satisfying $\tilde{f}^* \theta = \theta$, which is given by (4.5).
Then $\Gamma(\tilde{f})= \{ ((x,\xi),(y,\eta)) | 
\tilde{f}(x,\xi)=(y,\eta) \}$ is conical and Lagrangian according
to lemma 4.5.

Uniform families $a=\{a_m\}_{m \in M}$ of local symbols are already
defined but we consider here a slight generalization of (5.5) admitting
additional dependence of a second variable $y$, i. e. we require
\be
  |\partial^\alpha_\xi \partial^\beta_x \partial^\gamma_y a_m(x,y,\xi)|
  \le C_{\alpha,\beta,\gamma} (1+|\xi|)^{q-|\alpha|} ,
\ee
where $C_{\alpha,\beta,\gamma}$ is independent of $m$. We write 
${\cal U}S^q(B \times B \times {\bf R}^n \setminus 0)$ for all such
symbols $a=\{a_m\}_m$. Consider $\varphi=\{ \varphi_m \}_{m \in M}$
with the following properties. 
Each $\varphi_m : 
B \times B \times {\bf R}^n \setminus 0 \longrightarrow {\bf R}$ 
is a smooth map, positive homogeneous of degree one with respect to $\xi
\in {\bf R}^n \setminus 0 $, i. e. $\varphi_m(x,y,\tau\xi)=\tau \cdot
\varphi_m(x,y,\xi)$ and 
$d_{x\xi} \varphi_m, d_{y\xi} \varphi_m$ are $\neq 0$ on the canonical 
support of $a_m$ . Furthermore, the map
\bea 
  \{ (x,y,\xi) \in \en \mbox{canonical support of} \en a_m ~| ~
  d_\xi \varphi_m (x,y,\xi)=0 \}  \nonumber \\
  \longmapsto \{ (x,d_x \varphi_m(x,y,\xi), y , - 
  d_y \varphi_m(x,y,\xi)) ~|~ (x,y,\xi) \en \mbox{as above} \en \} \nonumber
\eea
is a diffeomorphism onto a conical submanifold     
$C_B \subset T^*B \setminus 0 \times T^*B \setminus 0$, where
$C_B$ corresponds to $C \subset T^*M \times T^*M$ under the
exponential map.

Such a family $\varphi=\{ \varphi_m \}_{m \in M}$ is called a {\em uniform
family of phase functions} associated to $a=\{a_m\}$ and we write 
${\cal U}Ph (a) (B \times B \times {\bf R}^n)$ for the space of all such
families. We say $A$ is a uniform Fourier integral operator of order
$q$, associated to the homogeneous canonical relation 
$C \subset T^*M \setminus 0 \times T^*M \setminus 0$, written as
$A \in {\cal U} F^q(M,C)$, if it satisfies the following conditions:

1. $A$ is a  continuous linear map $A: C^\infty_c(M) \longrightarrow
C^\infty(M)$,

2. If the family $\{ A_m \}_{m \in M}$ is defined by the commutative
diagram

\bi                                                   
                                                      
\[ C^\infty_c (B_\varepsilon(m))                      
   \quad                                              
   \begin{array}{c} A \\                              
   \Large\longrightarrow \\ {}                       
   \end{array}                                       
   \quad                                             
   C^\infty (B_\varepsilon(m)) \]                    
\vspace{-2ex}                                        
\[ \exp^*_m \Bigg\downarrow \hspace{4cm}             
   \Bigg\downarrow \exp^*_m \]                       
\vspace{-1ex}                                        
\[                                                   
   C^\infty_c(B)                                      
   \quad \quad                                        
   \begin{array}{c} A_m \\                            
   \Large\longrightarrow \\ {}                        
   \end{array}                                        
   \quad \quad                                        
   C^\infty(B) \]                                     
                                                      
\bi           
\no            
then there exist $a=\{ a_m \}_m \in {\cal U} S^q
(B \times B \times {\bf R}^n \setminus 0)$ such that
\be            
  A_m u(x) = (2\pi)^{-n} \int\limits_{{\bf R}^n} \int\limits_B
  e^{i\varphi_m(x,y,\xi)} a_m(x,y,\xi) u(y) dy d\xi, 
\ee             
$supp \en u \subset B$, and for $\psi :=(\exp^*_m)^{-1} \varphi_m$
the map
\be
  (x,y,\xi) \longmapsto (x, d_x \psi (x,y,\xi), y, - 
  d_y \psi (x,y,\xi))
\ee
is a diffeomorphism from the zero set of $d_\xi \psi$ onto a
submanifold of $C$.

Because  of our applications, we restrict ourselves to the case where
$C=\Gamma(\tilde{f})$, 
$f \in {\cal D}^{r+1}_{\theta,0} (S(T^*M))$, and we write simply 
${\cal U} F^q(f)$ for the corresponding class of uniform Fourier integral
operators. We set
\bea                                                                   
  {\cal U} F^{-\infty} (f) &:=& \bigcap\limits_q  {\cal U} F^q (f) , \\
  {\cal U} F^q (r+1) &:=& \bigcup\limits_{f \in 
  {\cal D}^{r+1}_{\theta,0} (S T^*M)}  
  {\cal U} F^q (f) , \\   
  \mbox{for} \en -k \le q ~,~~ {\cal U} F^{q,k} (f) &:=& {\cal U} F^q (f) /
  {\cal U} F^{-k-1} (f) , \\[1ex]
  {\cal U} F^{q,k} (r+1) &:=& \bigcup\limits_{f \in 
  {\cal D}^{r+1}_{\theta,0}}                      
  {\cal U} F^{q,k} (f) .                   
\eea                         
For $A_1 \in {\cal U} F^{q_1} (f_1)$ and $A_2 \in {\cal U} F^{q_2} (f_2)$
we have
\be
  A_1 \circ A_2 \in {\cal U} F^{q_1+q_2} (f_1 \circ f_2) .
\ee
Denote analogously to (5.18) for $-k \le q$
\be
  {\cal U} \Psi^{q,k} := {\cal U} \Psi^q / {\cal U} \Psi^{-k-1}
\ee
and 
${\cal D}^\infty_{\theta,0} = \lim\limits_\leftarrow 
{\cal D}^{r+1}_{\theta,0}$  , 
\bea
   {\cal U} F^{q,k} (\infty) &:=& \lim\limits_\leftarrow 
   {\cal U} F^{q,k} (r+1) \en = \en \bigcup\limits_{f \in 
   {\cal D}^{\infty}_{\theta,0}}                      
   {\cal U} F^{q,k} (f) . \\
   {\cal U} S^{q,k} &:=& {\cal U} S^q / {\cal U} S^{-k-1} .
\eea
Let $({\cal U} F^{0,k} (r+1))_*~$, $({\cal U} F^{0,k} (\infty))_*$ and
$({\cal U} \Psi^{0,k})_*$ denote the groups of invertible elements of 
${\cal U} F^{0,k} (r+1)$, ${\cal U} F^{0,k} (\infty)$ and 
${\cal U} \Psi^{0,k}$ respectively. It is clear from (5.20) that we must
choose $q=0$ to get invertibility inside one homogeneous constituent
of our graded structures.

\me
\no
{\bf Lemma 5.9.} {\it
Assume the hypothesizes of 4.8. The following  is an exact sequence of groups
\be
  I \longrightarrow ({\cal U} \Psi^{0,k})_*                
  \begin{array}{c} j \\[-2ex] \longrightarrow \\[-2ex] {} \end{array}   
  ({\cal U} F^{0,k}(\infty))_*     
  \begin{array}{c} \pi \\[-2ex] \longrightarrow \\[-2ex] {} \end{array}   
  {\cal D}^\infty_{\theta,0} \longrightarrow e ,
\ee 
where $j$ is the inclusion and $\pi[A]=f$ , where $A \in {\cal U}F^0(f)$,
$[A] = A + {\cal U}F^{-k-1}(f)$.
}

\me
\no
{\bf Proof.}
The injectivity of $j$, the surjectivity of $\pi$ and $im \, j \subseteq
ker \, \pi$
are clear. $im \, j = ker \, \pi$ follows from (5.26).
\hfill $\Box$

\me
\no
Generalizing the ideas of \cite{ARS1} and \cite{ARS2} our strategy is as
follows:

1. We want to construct a local section of $\pi$ in (5.24). For this we
need a chart in ${\cal D}^{r+1}_{\theta,0}$ at $id=e$.

2. This yields a chart and local section in (5.24), 
                                    
3. We endow $({\cal U} \Psi^{0,k})_*$ with the structure of an ILH Lie
group by forming Sobolev completions and taking the inverse limit.
                                                                  
4. We endow $({\cal U} F^{0,k}(\infty)_*$ with the structure of an ILH
Lie group, using these structures of $({\cal U} \Psi^{0,k})_*$ and 
${\cal D}^\infty_{\theta,0}$ in (5.24), the local section of $\pi$ and some
group theoretical theorems presented in section 7. 

\me
\no
This procedure is carried out in \cite{ARS1}, \cite{ARS2} for compact
manifolds, but in the case of open  manifolds the analysis is much harder.
We first start with the  construction of a local section of $\pi$ in (5.24).
This means the  existence of a neighborhood ${\cal U} (id) \subset 
{\cal D}^{r+1}_{\theta,0}$ and a (at least continuous) map 
$\sigma: {\cal U} \longrightarrow ({\cal U} F^{0,k}(r+1))_*$ such that
$\pi \circ \sigma = id_{\cal U}$. For doing this we construct global 
phase functions and present an explicit formula for $\sigma$.

We call a $C^2$ function $\varphi : T^*M \times M \rightarrow {\bf R}$
a {\em global phase function} for $e=id \in {\cal D}^{r+1}_{\theta,0} $
iff 

1. $d \varphi(T^*M \times M)$ is transversal to
$N_\pi := \{ \alpha \in T^*(T^*M \times M) | \alpha (v) = 0$ for all
$v \in ker \, \pi_* \} \subset T^*(T^*M \times M),$ where $ \pi: T^*(T^*M
\times  M)
\longrightarrow M \times M, \pi (\alpha_x,y)=(x,y)$

\no
and

2. $d \varphi (T^*M \setminus 0 \times M)_{N_\pi} := 
[ d \varphi (T^*M \setminus 0 \times M) \cap N_\pi ] / {\cal N}^\perp_\pi
= \Gamma (e) \subset (T^*M \setminus 0 \times T^*M \setminus 0,
\omega \ominus \omega )$, where ${\cal N}^\perp_\pi$ is the foliation
by isotropic submanifolds generated by the $\omega$--orthogonal
bundle $(TN_\pi)^\perp$ in $TT^*(T^*M \times M)$.

\me
\no
{\bf Example.}
For $M^n={\bf R}^n$ the function 
$\varphi: T^*{\bf R}^n \times {\bf R}^n \longrightarrow {\bf R}$
$\varphi(\xi_x,y)=<\xi,x-y>$, is a global phase function for $e$.
Consider $0<\delta<r_{inj}(M,g)$,
$\Omega_\delta = \{ (\xi_x,y) \in T^*M \times M | d(x,y)<\delta \}
\subset T^*M \times M$. In the definition above it is possible to 
replace $\varphi : T^*M \times M \rightarrow {\bf R}$ by
$\varphi |_{\Omega_\delta} : \Omega_\delta \subset T^*M \times M \rightarrow
{\bf R}$ if
$\Omega_\delta$ contains the whole fibers of $T^*M$.

\me
\no
{\bf Lemma 5.10.} {\it
The function $\varphi_0 : \Omega_\delta \rightarrow {\bf R}$ defined by 
$\varphi_0(\xi_x,y)=\xi_x (\exp^{-1}_x (y))$, is a global phase
function for $e =id \in {\cal D}^{r+1}_{\theta,0}$ on $\Omega_\delta$.
}

\me
\no 
We refer to \cite{ARS1}, p. 541,42 for the proof which is local in character
and does not depend on the compactness or openness of $M$.
\hfill $\Box$
   
\me
For $\delta < r_{inj} (M,g)$ we can define global
symbols $a(x,\xi) \in {\cal U} S^q (\Omega)$ on 
$\Omega = \Omega_\delta$                           
by requiring 
\be
  |\nabla^\alpha_\xi \nabla^\beta_x a(x,\xi)| \le C_{\alpha,\beta}
  (1+|\xi|^{q-\alpha}
\ee
and assuming an asymptotic expansion of the type (5.10), (5.11).

Let $\chi(x,y)$ be a bump function on $M \times M$ such that 
$supp \, \chi \subset {\cal U}_\delta (diagonal) = 
{\cal U}_\delta (\bigtriangleup)$ 
and $\chi \equiv 1$ on
a neighborhood of the diagonal.

\me
\no
{\bf Proposition 5.11.} {\it 
Let $a(x,\xi) \in {\cal U} S^q (\Omega)$ be a global classical 
symbol of order $q$. Then 
\be
   u(x) \longmapsto A u(x) := (2 \pi)^{-n} 
   \int\limits_{T^*_xM} \int\limits_{B_\delta(x)}
   \chi(x,y) e^{i\varphi_0(\xi_x,y)} a(x,\xi)
   |det \exp_{x*}| dy d\xi
\ee
is a classical pseudodifferential operator of order $q$ in the former
sense, i. e. $A \in {\cal U}\Psi^q(m)$.
}

\me
\no
{\bf Proof.} The fact that $A$ is a $\Psi$DO is well known and follows
from the famous Kuranishi principle. The uniform boundedness
(5.5) of the family $\{ a_m \}_m$ follows from 5.10 and 5.25.
\hfill $\Box$

\me
Now we return to our first task to define a chart in 
${\cal D}^{r+1}_{\theta,0}$. The most simple idea would be to endow
${\cal D}^{r+1}_0 \times {\cal F}^{r+1}_0$ with a Riemannian metric 
$G$, to take the induced metric $G_\theta$ on the submanifold
${\cal D}^{r+1}_{\theta,0}$ and then apply the Riemannian exponential
map of $G$ to a sufficiently small ball in 
$T_{(id,1)} {\cal D}^{r+1}_{\theta,0}$. But this will not work since
${\cal D}^{r+1}_{\theta,0}$ is definitely not a geodesic submanifold, 
at least for any reasonable metric on 
${\cal D}^{r+1}_0 \times {\cal F}^{r+1}_0$.
We could take the Riemannian exponential of $G_\theta$, but we don't
know it, i. e. we can't calculate or estimate this. Therefore we will
construct a chart centered at 
$id: T^*M \setminus 0 \rightarrow T^*M \setminus 0$ by another method.
The framework of our approach here is already carried out for 
compact manifolds in 4.2, 4.3, 4.4
of \cite{ARS1},  but the case $M^n$ open has its own features. 
It requires additional estimates at $\infty$.

Assume $(W,g)$ satisfies the conditions $(I)$ and
$(B_{\infty})$. 
Let $y\in W, ~\delta < r_{inj}$
and $(U^{\delta}(y), y^1,...,y^n)$ a normal chart.
Let
$v=v^1\frac{\partial}{\partial y^1}+ \cdots+
v^n\frac{\partial}{\partial y^n}$ be a locally defined vector field, 
$| v|=(v^iv_i)^{\frac{1}{2}}$ its Riemannian norm and 
$|v|_{fl}=(\sum_{i=1}^n(v^i)^2)^{\frac{1}{2}}$ its euclidean flat 
pointwise norm respectively.
We want to compare $|v|$ and $|v|_{fl}$.
More generally, consider a locally defined tensor field 
$t=(t^{i_1 \cdots i_k}_{j_1\cdots j_l})$.
Denote as before by $|t|$ its Riemannian norm and by 
$|t|_{fl}= (\sum (t^{i_1 \cdots i_k}_{j_1\cdots j_l})^2)^{\frac{1}{2}}$ its
flat euclidean norm.
Moreover , denote by $|t|_{U,r+1,fl}$ the $(r+1)$ - Sobolev norm based on 
Riemannian covariant derivatives, the Riemannian volume element and the
flat pointwise norm.

\me
\no
{\bf Lemma 5.12.} a)
{\it There exists constants $c_1 ,c_2 >0$ such that for any
$(k,l)$-tensor $t$ 
\be \label{a}
c_1|\nabla^it|_{x,fl}\leq|\nabla^it|_x\leq c_2|\nabla^it|_{x,fl}
\ee
for all $x\in U_{\delta}, 0\leq i \leq r+1, , c_1, c_2 $ independent of
$y$, depending only on $k,l, r$ . }\\
b)
{\it If $t\in\Omega^{0,2,r+1}(T^k_lW)$ and ${\cal U}=\{(U_{\lambda}, 
\phi_{\lambda})\}_{\lambda}$ is an uniform cover of $W$ by normal charts
of radius $\delta$ then there exists constants $d_1 , d_2 >0$ such that
\be
d_1 \sum_{\lambda} |t|_{U_{\lambda}, r+1,fl} \leq|t|_{r+1} \leq 
d_2 \sum_{\lambda} |t|_{U_{\lambda}, r+1,fl} .
\ee
 }

\me
\no
{\bf Proof.}
a) In $U_{\delta}(x)$ holds
\be
k_1\delta_{ij}\xi^i\xi^j \leq g_{ij}\xi^i\xi^j \leq k_2
\delta_{ij}\xi^i\xi^J
\ee
\be
k_3\delta^{ij}\xi_i\xi_j \leq g^{ij}\xi_i\xi_j \leq
k_4\delta^{ij}\xi_i\xi_j .
\ee

b) This follows from a) using a partition of unity that is bounded up to
order $r+1$ , the module structure theorem and the fact that the cover is
uniformly locally finite.

\hfill $\Box$

\me

Now we generalize our notion of a global phase function $\varphi_H$ for
$e=id\in  {\cal D}_{\theta ,0}^{r+1} (T^*M \setminus 0)$ from above to
other 
$f\in {\cal D}_{\theta ,0}^{r+1} (T^*M \setminus 0)$.
Consider again the projection $\pi : T^*M \times M \longrightarrow M ,
~\pi(\alpha_x ,y)=(x,y)$, the conormals $N_{\pi} \subset T^*(T^*M\times
M)$ of this submersion and the corresponding foliation ${\cal
N}^\perp_{\pi}$ of
$N_{\pi}$.
A function $\varphi : T^*M \times M \longrightarrow {\bf R}$ is called a 
{\em global phase function} for a diffeomorphism $f$  of $T^*M \setminus 0$
iff:

1. $d\varphi (T^*M\times M) \subset T^*(T^*M\times M)$ is transversal to
$N_{\pi}$, and

2. $\Gamma (f) = d\varphi (T^*M\times M)_{N_{\pi}} := [(d\varphi
(T^*M\times M)\cap N_{\pi})/{\cal N}^\perp_{\pi}]'$,  \\ where 
$C':=\{(x,y,\xi,-\eta )~|~(x,y,\xi \eta ) \in C\}$.

For our purposes we can restrict ourselves to maps 
$\varphi :\Omega_{\delta} \subset T^*M \times M \longrightarrow {\bf
R}$.\\
In local charts we have the following representations:\\ 
$N_{\pi}=\{(x,\alpha ,\vartheta, 0,y,\eta )\}$ \\
${\cal N}^\perp_\pi =\{ (x, {\bf R}^n,\vartheta ,0,y,\eta )\}$ = one leaf
\\
$d\varphi |_{(\alpha_x ,y)}= (x,\alpha ,d_x\varphi , d_\alpha \varphi , y,
d_y \varphi) $,\\
hence we have locally

\bea \label{A}
d\varphi (\Omega )_{N_\pi}&=&(d\varphi (\Omega )\cap N_\pi )/{\cal
N}^\perp_\pi \nonumber \\
 &=&\{(x,\alpha ,d_x\varphi ,d_\alpha \varphi =0,y,d_y\varphi )~|~(\alpha
,y)\in\Omega\} \nonumber\\
 &=&\{(x,y,d_x\varphi ,d_y \varphi )~|~d_\alpha \varphi =0\} .
\eea

We prescribe a $\varphi$ and want to construct the corresponding $f$. 
For $(x,y)\in U_ \delta (\Delta )\subset M\times M $ denote
$v(x,y)=\exp_x^{-1}(y)\in  T_xM$ and as above
$\varphi_0(\alpha_x,y)=<\alpha ,v(x,y)>$.

\me
\no
{\bf Lemma 5.13.} {\it 
Let $h \in \Omega^{0,2,r+2}(S)$, $|h|_{r+2}<\varepsilon$, 
$\varepsilon$ sufficiently small and $H$ be the extension homogeneous of
degree one to $T^*M \setminus 0$. 
Define $\varphi_H : \Omega \subset T^*M
\setminus 0 \times M
\rightarrow {\bf R}$ by
\[ \varphi_H(\alpha_x,y) := \varphi_0(\alpha_x,y) + H(\alpha_x) =
<\alpha ,v(x,y)>+H(x,\alpha ). \]
Then there exists an $f \in {\cal D}^{r+1}_{\theta,0} (T^*M \setminus 0)$ 
such that $\varphi_H$ is a global phase function for $\Gamma(f)$
i. e. $d \varphi_H((T^*M \setminus 0) \times M)_{N_\pi} = \Gamma(f)$.
}

\me
\no
{\bf Proof.} Let $y\in M$ and $e_1,...,e_n$ be an orthonormal base of
$T_yM$. Let $\delta <R_{inj}$ and 
$(U^{\delta}(y), y^1,...,y^n)$ the corresponding  normal chart.
Then $v(x,y)\in T_xM$ and $\alpha\in T^*_xM$ can be written uniquely as 
$v=v^1\frac{\partial }{\partial y ^1}|_x + \cdots
+v^n\frac{\partial }{\partial y ^n}|_x$ and $\alpha=\alpha_1dy^1+\cdots +
\alpha_ndy^n$ respectively.
The transversality property is substantially a local property and has
been established in \cite{ARS1}.
We must show that $[d\varphi_H (\Omega )_{N_\pi}]'=\Gamma (f)$ for some 
$f\in {\cal D}_{\theta ,0}^{r+1} (T^*M \setminus 0)$, if
$|h|_{r+2}<\varepsilon$, for $\varepsilon$ sufficiently small.
The Lagrangian submanifold in $T^*M \setminus 0 \times T^*M \setminus 0$
generated by $\varphi_H$ can, according to (\ref{A}), locally be written as
\[
\{(x,y,d_x\varphi_H,d_y\varphi_H~|~d_\alpha\varphi_H=0\} .
\]
We assume that $(y,\eta )\in S(T^*M)$ be given and we have to solve for
$(x,\alpha )$, i.e we have to solve the equations
\be \label{B}
d_y\varphi_H=d_y<\alpha ,v>=<\alpha ,d_yv>=-\eta
\ee
\be \label{C}
d_\alpha\varphi_H=0 .
\ee
(\ref{B}) yields
$$\left( \begin{array}{c}
\alpha_1\\
.\\
.\\
.\\
\alpha_n
\end{array}\right) =-
\left( \begin{array}{ccc}
\frac{\partial v^1}{\partial y^1}& \cdots & \frac{\partial v^n}{\partial y^1} \\
.& & . \\
.& &. \\
.& &. \\
\frac{\partial v^1}{\partial y^n}& \cdots & \frac{\partial v^n}{\partial y^n}
\end{array} \right)_{| y}^{-1}
\left( \begin{array}{c}
y_1\\
.\\
.\\
.\\
y_n
\end{array}\right) ~,
$$
where $\eta=\eta_1dy^1 |_y + \cdots +\eta_ndy^n |_y$ . 
The matrix $\left( \cdots \right)^{-1}$ equals to 
$\frac{\partial(y^1,...,y^n)}{\partial (v^1,...,v^n)}$, i.e we obtain
\be\label{D}
\alpha_i=-\frac{\partial y^j}{\partial v^i}\eta^j~,~~~
\alpha=-\frac{\partial y^j}{\partial v^i}\eta^j dy^i |_x~.
\ee
(\ref{D}) can be reformulated as follows: 
The map $v:(x,y) \mapsto \exp_x^{-1} (y)$ maps $U_\delta (\Delta )$ into $TM$.
Fixing (the unknown) $x$ for a moment, we obtain a map 
$\exp_x^{-1} ( . ) : U_\delta(x) \longrightarrow T_xM$, its $y$ - differential
maps $T_yU_\delta (x) \longrightarrow T_{v(x,y)}  T_xM $ and 
$(v(x, .)_{*,y})^{-1}=(v(x,.)^{-1})_{*,y} = (\exp_x)_{*,y} :T_{v(x,y)}T_xM
\longrightarrow T_y U_\delta$, i.e we can write
\be \label{E}
\alpha=-\eta\circ (\exp_x)_{*,y}
\ee
The equation (\ref{C})  then becomes
$$
d_\alpha\varphi_H =d_\alpha <\alpha , v(x,y)>+d_\alpha H=0  ,
$$
\be \label{Ea}
v^1d\alpha^1+\cdots v^nd\alpha^n=-(\frac{\partial H}{\partial\alpha_1}d\alpha_1
+\cdots \frac{\partial H}{\partial\alpha_n}d\alpha_n )|_{-\eta v_{*,y}^{-1}}~ ,
\ee
which we have to solve for $x$.
Equation (\ref{Ea}) as an equation between locally defined
vector fields is equivalent to $v^i=-\frac{\partial H}{\partial \alpha_i}|_{-\eta
v_{*,y}^{-1}}~~,1\leq i\leq n$, which is equivalent to 
\be \label{Eb}
v^1\frac{\partial}{\partial y^1}|_x +\cdots +v^n\frac{\partial }{\partial y^n}|_x
= -(\frac{\partial H}{\partial\alpha_1}\frac{\partial }{\partial y^1}|_x+\cdots +
\frac{\partial H}{\partial\alpha_n}\frac{\partial }{\partial y^n}|_x)|_{-\eta
v_{*,y}^{-1}}~,
\ee
which we have to solve for $x$.
There are several way to solve (\ref{Eb}).
Write the right hand side of (\ref{Eb}) as 
$-<\frac{\partial H}{\partial\alpha } ~,~ \frac{\partial }{\partial y}|_x >$ , we
can assume $c_2 |<\frac{\partial H}{\partial\alpha } ~,~ \frac{\partial }{\partial
y} >|_{fl} =\delta_1 <\delta$, (cf.(\ref{a})).
consider the map 
$F: [0,1]\times U_\delta (y) \longrightarrow V_\delta (0_{TU_\delta })~,~F(t,x) =
\exp_x^{-1}(y) +t<\frac{\partial H}{\partial\alpha } ~,~ \frac{\partial }{\partial
y}|_x >$ and set 
$L:=\{t\in [0,1]~|~there~ exists~ a ~unique~ x_t~ such~ that~ F(t,x_t)=0\}$. 
Then the following hold:\\
1) $L\neq\emptyset$, since $0\in L$ with $x_0=y $.
2) There exists a $\lambda >0$ such that $[0,\lambda [ \subset L$. This follows from
the fact that $d_xF|_{(0,y)}$ is invertible.
3) $L$ is open, which follows analogous to 2).
4) $L$ is closed: assume $t_1<t_2<\cdots\rightarrow t^*$, this implies
$x_{t_\nu} \in U_{t^*\delta}$, there exists a convergent subsequence $x_{t_{\nu_i}}
\rightarrow x^*$ and $F(t^*,x^*)=0$. 
Therefore $L=[0,1]$ and there exists a unique
$x\in U_\delta (y)$ such that $\exp_x^{-1}(y)=-<\frac{\partial H}{\partial\alpha }
~,~ \frac{\partial }{\partial y}|_x >|_{\alpha ={-\eta v_{*,y}^{-1}}}~.$\\
What remains to assure is that $c_2 |<\frac{\partial H}{\partial\alpha }
~,~ \frac{\partial }{\partial y} >|_{fl}<\delta$.
This can been seen as follows: For a covector $V\in T^*T^*M$ we have 
$|V|=g(KV,KV)+g(V_h,V_h)$ where $K:T^*T^*M\longrightarrow T^*_{vert}T^*M$ is the
connection projection.
Any $V\in T^*T^*U_\delta$ can be written as 
$V=A_1dy^1+\cdots +A_ndy^n+B_1d\alpha^1+\cdots +B_nd\alpha^n$.
If $V=B_1dy^1+\cdots +B_ndy^n$, then it is purely vertical and
$g_{Sasaki} (V,V)=g(KV,KV)=g(V,V)$. 
This can be calculated from the $g^{ij}(x)$ and the $B_i's$ i.e. we can apply
(\ref{a}).
Using uniform boundedness, properties of $\exp$ (cf.\cite{ES}) we derive from
(\ref{a}) and (\ref{D}),(\ref{E})
\be\label{F}
\frac{1}{2}=\frac{1}{2}|\eta |_y \leq |\alpha |_x \leq\frac{3}{2}|
\eta|_y =\frac{3}{2}
\ee
if we choose $\delta<r_{inj}$ small enough.
If we assume this has been done, $\alpha$ lies in a $\frac{1}{2}$- neighborhood of
$S(T^*M)$.
By assumption $H|_S=h\in \Omega^{0,2,r+2}(S)$,
\[
{}^b|h|\leq C_0|h|_{r+2}
\]
\be \label{G}
{}^b|d_\alpha H|_{U_\frac{1}{2}(S)}\leq
~{}^b|dH|_{U_\frac{1}{2}(S)}\leq\frac{1}{2}~{}^b|d_Sh|+{}^b|h|\leq 
\frac{1}{2}C_1|d_Sh|_{r+1}+C_0|h|_{r+2}\leq C_2|h|_{r+2} ~,
\ee
hence
\bea \label{H}
|v(x,y)|_x=|\exp_x^{-1}(y)|\leq c_2|<\frac{\partial H}{\partial\alpha }
~,~ \frac{\partial }{\partial y} >|_{fl}\leq c_2c_1|d_\alpha H|_y \nonumber\\
\leq 
c_2c_1~ {}^b|d_\alpha H|_{U_\frac{1}{2}(S)} \leq c_2c_1
~{}^b|dH|_{U_\frac{1}{2}(S)}
\leq c_2c_1 C_2|h|_{r+2}=C_3|h|_{r+2} .
\eea
Choosing $|h|_{r+2} \leq \frac{\delta}{C_3}$ we obtain a unique $x$ solving
equation (\ref{C}).
We get an orthonormal base at $x$ by parallel translation of
$e_1,...,e_n\in T_yU_\delta$ from $y$ to $x$ along the unique connecting geodesic.
This yields a normal chard $(U_\delta (x), x^1,...,x^n)$, $y\in U_\delta (x)$.
We now insert $x$ into $<\alpha , d_xv(x,y)>+d_xH(x,\alpha )$ and obtain
\bea \label{I}
\xi (x,y)&=& <\alpha , d_xv(x,y)>+d_xH(x,\alpha ) \nonumber\\
 & =& -(\frac{\partial y^j}{\partial v^1}\eta^j d_xv^1+\cdots + 
\frac{\partial y^j}{\partial v^n}\eta^j d_xv^n )+d_xH \nonumber\\
 &=& -\frac{\partial v^k}{\partial x^l}|_x\frac{\partial y^j}{\partial 
v^k}|_{v^{-1}(y)} \eta_j dx^l+\frac{\partial H}{\partial x^l}|_xdx^l ~.
\eea
The coordinate free description of $<\alpha , d_xv(x,y)>$ is given by
\bea \label{J}
<\alpha , d_xv(x,y)>&=&-\beta\circ [(\exp_x)_{*y}:T_{v(x,y)}T_xM\rightarrow 
T_y U_\delta (y)] \circ\nonumber\\
&\circ& [(\exp^{-1}_{(.)}(y)_{*y}:T_xU_\delta (x)\rightarrow T_{v(x,y)}T_xM~.
\eea
If $0\neq\eta\not\in S(T^*M)$, then $\alpha =-\eta\circ v_{*,y}^{-1}$ remains , but
to determine $x$ we apply the procedure above with $\eta ' =\eta /|\eta |$ and
$\alpha ' =\eta ' \circ v_{*,y}^{-1}$ whereas $\xi$ is gain defined by (\ref{I}).
Thus we finally get a map $f^{-1}(y,\eta )=(x(y,\eta ), \xi (y,\eta ))$.
This map is well defined .
The determining equations and its solutions can be described coordinate free as
shown by (\ref{E}) and (\ref{J}).
Using the geodesic convexity of normal charts and the implicit function theorem,
it is easy to derive that $f^{-1}$ is $1-1$, onto and of class $C^1$ together with
its inverse.
$f$ and $f^{-1}$ are $C^1$ - diffeomorphisms of $T^*M\setminus 0$.
The fact that $d\varphi _H (T^*M\setminus 0\times M)$ is a conic Lagrangian
submanifold implies that $f^*\theta =\theta$, hence $f\in C^1{\cal D}_\theta
(T^*M\setminus 0)$.
To assure that $f\in {\cal D}^{r+1}_{\theta ,0}(T*M\setminus 0)$ we have to show
that the map $f^{-1}:S(T^*M)\ni\eta_y\longrightarrow \xi_x /|\xi_x| \in S(T^*M)$
belongs to ${\cal D}^{r+1}_{\theta , 0}(S(T^*M))$.
To assure that this new $f^{-1}\in {\cal D}^{r+1}_\theta (S(T^*M))$ it would be
sufficient that $dist((x, \xi_x /|\xi_x|),(y,\eta ))<r_{inj}S(T^*M)$ and that 
the corresponding vector field on $S(T^*M)$ would be of class $r+1$.
Hence we have to estimate form above $dist(\xi /|\xi | , \eta )$.
For this we need an arc from $\xi$ to $\eta$.
First consider the parallel translation of $\xi /|\xi |$ to the fiber  $S(T^*_yM)$
along the geodesic $\exp (sv(x,y))$, $0\leq s\leq 1~$,$~s\mapsto P_s\xi/|\xi |$, 
$~P_1\xi /|\xi |\in S(T^*_yM)$. 
This is a horizontal geodesic in $S(T^*M)$ covering the geodesic $\exp (sv)$.
We obtain $length\{P_s\xi /|\xi |\}=length\{\exp (sv)\}=|v|$ ,$~dist(\xi ,
P_1\xi/|\xi |) = dist(x,y)=|v|$.
But $P_1 \xi /|\xi |$ and $\eta$ lie in the fiber $S(T^*_yM)$ which is an
euclidean sphere.
The distance of two points is the length of the shortest geodesic between them.
But this distance can be estimated from above by the distance in $T^*_yM$
multiplied by a factor $C_4\approx 1.8$, i.e. 
$dist(P_1\xi /|\xi | , \eta )\leq C_4|P_1\xi / |\xi |  - \eta |$.
Now 
\bea 
|P\xi /|\xi | - \eta |= \left| \frac{P\xi |\eta | - \eta |\xi |}{|\xi |}\right|
= \left|\frac {P\xi |\eta | -P\xi |\xi | + P\xi |\xi | -\eta |\xi |}{|\xi |}\right|
\nonumber \\
\leq \left|\frac{P\xi (|\eta |-|\xi |)}{\xi |}\right|+|P\xi - \eta | 
=
|~|\eta |-|\xi |~|+|P\xi -\eta |\leq 2|P\xi - \eta | ~.\nonumber
\eea
If $\xi =\xi_idx^i|_x $ and $P\xi = (\xi_i+ \Delta \xi_i)dx^i|_y$ 
then $\xi_i +\Delta \xi_i$ is the solution $\xi_i(|v|)$ of the equation 
$$\frac{d\xi_i(t)}{dt}-\Gamma^k_{ij}\xi_k(t)\frac{dx^j(t)}{dt}=0 ,$$
where $x^j(t)$ are the coordinates of the geodesic
$\exp(\frac{t}{|v|}v)~,~\xi_i(0)=\xi_i$. 
It is a simple and well known fact that $\Delta \xi |$ can be estimated as
\be\label{K}
|\Delta \xi |\leq C_5 \cdot\Gamma \cdot |v|\cdot |\xi| =C_6 ~dist(x,y)\cdot |\xi|
\ee
Hence
\[
|P\xi - \eta |_y = |\xi_idx^i_{|y} + \Delta \xi_idx^i_{|y} - \eta | 
 \leq   |\xi_idx^i_{|y} - \eta |+C_6~dist(x,y)\cdot |\xi | ~.
\]
We insert $\xi_ldx^l$ from (\ref{I}) and $\eta =\eta_idy^i|_y =\frac{\partial
y^j}{\partial x^l}\eta_jdx^l|_y$ and obtain
\bea\label{L}
|\xi_ldx^l_{|y} -\eta |&\leq & \left|-\frac{\partial v^k}
{\partial x^l}|_x \frac{\partial y^j}{\partial v^k}|_{v^{-1}(y) }~\eta_j -
\frac{\partial y^j}{\partial x^l}|_y \eta_jdx^l \right|_y 
+\left| \frac{\partial H}{\partial x^l}|_xdx^l|_y\right|_y \nonumber \\
&\leq & \left|-\frac{\partial v^k}
{\partial x^l}|_x \frac{\partial y^j}{\partial v^k}|_{v^{-1}(y) }~\eta_j -
\frac{\partial y^j}{\partial x^l}|_y \eta_jdx^l \right|_{fl} 
+\left| \frac{\partial H}{\partial x^l}|_xdx^l|_y\right|_y~.
\eea
Next we use the following uniform expansions on $M$:
\[
\frac{\partial (y^1,...,y^n)}{\partial (v^1,...,v^n)}|_{v^{-1}(y)} =
\left(\begin{array}{ccc}
1& &  0\\
 & \ddots & \\
0&  & 1 
\end{array}\right)
+O(dist(x,y))+o(dist(x,y)) ,
\]
\[
\frac{\partial (v^1,...,v^n)}{\partial (x^1,...,x^n)}|_x =
-\left(\begin{array}{ccc}
1& &  0\\
 & \ddots & \\
0&  & 1 
\end{array}\right)
+O(dist(x,y))+o(dist(x,y)) ,
\]
\[
\frac{\partial (y^1,...,y^n)}{\partial (x^1,...,x^n)}|_{x(y)} =
\left(\begin{array}{ccc}
1& &  0\\
 & \ddots & \\
0&  & 1 
\end{array}\right)
+O(dist(x,y))+o(dist(x,y)) .
\]
Taking this into account and $|\eta |_y=1$ we can conclude that  the right hand side
of (\ref{L}) can be estimated as
\[
(\ref{L}) \leq C_7|\eta |~dist(x,y)+\frac{3}{2}|d_xH(x, \alpha /|\alpha |)| \leq
C_7 dist(x,y)+\frac{3}{2}~{}^b|d_Sh|\leq C_7~dist(x,y)+C_8|h|_{r+2}~.
\]
There remains to estimate $|\xi |$ itself .
\bea
|\xi |_x &\leq  & |<\alpha , d_x v>|_x+|d_xH(x, \alpha )|_x \leq c_2~|<\alpha , d_x
v>|_{fl} +C_8|h|_{r+2} \nonumber\\
 & \leq & c_2 (|\eta |+C_9 |\eta | dist(x,y))+C_8|h|_{r+2} =
c_2+C_{10}dist(x,y)+C_8|h|_{r+2}~. \nonumber
\eea

We now reached our final estimate :
\bea\label{M}
dist (\xi , \eta )&\leq& dist (x,y)+2C_4[C_6|\xi
|dist(x,y)+C_7dist(x,y)+C_8|h|_{r+2}] \nonumber \\
 &\leq&
dist(x,y)+2C_4[C_6(c_2+C+{10}dist(x,y)+C_8|h|_{r+2})dist(x,y) +\nonumber \\ 
& + &C_7dist(x,y)+C_8|h|_{r+2}]\nonumber\\
& \leq &
C_3|h|_{r+2}+2C_4[C_6(c_2+C_{10}C_3|h|_{r+2}+C_8|h|_{r+2})C_3|h|_{r+2} + \nonumber\\
&+&C_7C_3|h|_{r+2}+C_8|h|_{r+2}] \nonumber\\ 
& \leq & C_{11}|h|_{r+2}+C_{12}|h|^2_{r+2}\leq C_{13}|h|_{r+2}~~,~if~|h|_{r+2}\leq 1~.
\eea
We conclude, that if 
\be \label{N}
\varepsilon_1 = \min \left\{
1,\frac{\delta}{C_3},\frac{r_{inj}(S)}{C_{13}}\right\}~, and
~~|h|_{r+2}<\varepsilon_1
\ee
then there exists a unique vector field $X(h)$ such that
\be
f=\exp X(h)~.
\ee
To obtain $f\in {\cal D}_0^{r+1}(S)$ we must in addition show that
$|X(h)|_{r+1}<\infty$.
This is a rather long and technical estimate,
but the established inequalities carry over in a quite natural manner step by step to
pointwise norms of derivatives (here we use repeatedly (\ref{a})) and finally to Sobolev
norms, applying the module structure theorem. We omit the details here.
We established that $f\in C^1{\cal D}_\theta(T^*M\setminus 0)$ and 
$f|_S =\exp X(h)\in C^1{\cal D}_0(S)$.\\
If we define $\beta$ by (4.6) we obtain $(f|_S=\exp X(h)~,~\beta )\in 
{\cal D}^{r+1}_{\theta ,0} (S)$, i.e the extended 
$f\in {\cal D}^{r+1}_{\theta ,0}(T^*M\setminus 0)$. 
\hfill $\Box$

\me
Now we sharpen our considerations by proving 
 
\me
\no
{\bf Lemma 5.14.} 
{\it The map 
$h \longmapsto H \longmapsto f$ given by Lemma 5.13 is a bijection from a
neighborhood ${\cal V}={\cal V}_\varepsilon(0) \subset
\Omega^{0,2,r+2}(S)$ onto a neighborhood 
${\cal U} = {\cal U}(id) \subset {\cal D}^{r+1}_{\theta,0} 
(T^*M \setminus 0)$ of $~id \in {\cal D}^{r+1}_{\theta,0}(T^*M \setminus 0)$.
The inverse mapping is given by
\be
  f \longmapsto H=H_f, \en
  H(\alpha_x) = -\langle \alpha_x, \exp^{-1}_x (\pi (f^{-1} (\alpha_x)))
\rangle. 
\ee
}

\me
\no
{\bf Proof.} First we recall the local topology of 
${\cal D}^{r+1}_{\theta,0}(S)$, giving a neighborhood basis for
$id \in {\cal D}^{r+1}_{\theta,0}(S)$. As established in \cite{Ei3},
a neighborhood basis for $id \in {\cal D}^{r+1}_0(S)$ is given
by ${\cal U}= \{ {\cal U}_\frac{1}{\nu} (id) \}_{\frac{1}{\nu} \le
r_{inj}(S)}$, where
\be
   {\cal U}_{\frac{1}{\nu}} (id) =  \left\{ \exp X ~| \en 
   |X|_{r+1} < \frac{1}{\nu} \right\} , \en
   (\exp X)(x) = \exp_x X(x).
\ee
Then the corresponding basis in ${\cal D}^{r+1}_{\theta,0}(S)$ is given 
by ${\cal U}_{\frac{1}{\nu},\theta}(id) = 
\left\{ {\cal U}_{\frac{1}{\nu},\theta}(id) \right\}_
{\frac{1}{\nu} \le r_{inj}(S)}$,
where
\be
   {\cal U}_{\frac{1}{\nu},\theta}(id) = {\cal U}_{\frac{1}{\nu}}(id) 
   \cap {\cal D}^{r+1}_{\theta,0} (S) .
\ee

\no
{\bf Remark.} {\it
The $Xs$ in (5.49) must not necessarily satisfy $L_X \theta = 0$ or
$L_X \theta = u \cdot \theta$. We follow in (5.48), (5.49) the
Riemannian exponential, not the integral curves of $X$ with
$L_X \theta = u \cdot \theta$. } \hspace*{2cm} {} \hfill $\Box$ 

\me
Considering (5.48), we see that $H$ is only well defined if 
$\exp^{-1}_z(\pi f^{-1}(\xi_z))$ is well defined, i. e. if 
\be
  dist (z, \pi f^{-1}(\xi_z)) < r_{inj} (M) .
\ee
Clearly, $r_{inj}(S,g_S) \le r_{inj}(M,g)$, since horizontal
geodesics project isometrically to geodesics.

For $\varepsilon_2$ sufficiently small and 
$|X|_{r+1}<\varepsilon_2$, we have
\be
   dist(\xi_z, (\exp X)^{-1} (\xi_z)) < r_{inj}(S) \le 
   r_{inj} (M) ,
\ee
i. e. for $f \in {\cal U}_{\varepsilon, \theta}(id)$, (5.51) is 
satisfied and 
\[ H(\xi_z) = - \langle \xi_z, \exp^{-1}_z (\pi f^{-1}(\xi_z)) \rangle \]
is well defined. Moreover, the right hand side is homogeneous of 
degree one. We must still assure that 
$H|_S \in \Omega^{0,2,r+2}(S)$. The derivatives of 
$H|_S$ lead to the derivatives of certain Jacobi fields. Their
pointwise norms can be estimated by polynomials on the pointwise norms
of the derivatives of $X$ which are additionally square integrable. 
This has been performed in \cite{Ei3}. At the end we get 
\be
  H |_S \in \Omega^{0,2,r+2} (S) .
\ee
The gain of one Sobolev index comes from the fact $\exp X$ contains
already one integration. Moreover given any $\varepsilon_3>0$, there
exists $\varepsilon_4>0$ such that
\be
  |X|_{r+1}<\varepsilon_4 \en \mbox{implies} \en 
  |H|_S|_{r+2} < \varepsilon_3 .
\ee
The key for proving (5.54) are the Jacobi field constructions, their
estimates and the module structure theorem.

Setting $\varepsilon_3=\varepsilon_1$ from (5.46) and choosing
$\varepsilon_5=\min \{ \varepsilon_2, \varepsilon_4 \}$, we obtain
\be          
   H=H_f \en \mbox{is for} \en f \in {\cal U}_{\varepsilon_5, \theta} (id)
   \en \mbox{well defined,} \en H|_S \in \Omega^{0,2,r+2}(S)
   \en \mbox{and} \en |H|_S|_{r+2} \le \varepsilon_1 .
\ee          
(5.55) now permits to construct the sequence of maps
\be         
   f \longmapsto H=H_f \longmapsto \varphi_{H_f} 
   \longmapsto f_{\varphi_{H_f}}
\ee                                 
In \cite{ARS1} it is proved that (5.56) equals to $id$.
If we set $\varepsilon = \min \{ \varepsilon_1, \varepsilon_5 \}$,
the sequence 
\be        
   H \longmapsto \varphi_H      
   \longmapsto f = f_{\varphi_H} \longmapsto H_f
\ee                                 
is also well defined and equals to $id$ according to \cite{ARS1}. 
We omit the proof that the constructed 1--1 mapping
$\Psi: f \longmapsto H_f \longmapsto h_f = H_f|_S$
is of class $C^{k-r}$.              
\hfill $\Box$                       
                                    
\me
We summarize our result in

\me
\no
{\bf Theorem 5.15.} {\it
$({\cal U}={\cal U}_{\varepsilon,\theta}(id), \Psi, 
\Omega^{0,2,r+2}(S))$ is a $C^{k-r}$ chart at 
$id \in {\cal D}^{r+1}_{\theta,0}(T^*M \setminus 0)$, where
\be
  \Psi(f)(\alpha_x) = - \langle \alpha_x, \exp^{-1}_x (\pi (f^{-1}
(\alpha_x)))
   \rangle
\ee
for all $\alpha_x \in T^*M \setminus 0$.
}
\hfill $\Box$
                              
\me                           
Now we are in a position to construct a local section $\sigma : 
 {\cal U}\subset{\cal D}^\infty_{\theta,0} = \lim\limits_{\leftarrow}
{\cal D}^{r+1}_{\theta,0} \longrightarrow ({\cal U} F^{o,k}(\infty))$.
We assume the condition $(B_\infty)$. For $f \in {\cal U} \cap 
{\cal D}^\infty_{\theta,0}$ and $A \in {\cal U} F^{o,k}(f)$ and
$a(x,\xi)$ a representative of the classical symbol of $A$ we
now can write in analogy to (5.26)       
\be                              
  Au(x)=(2\pi)^{-n} \int\limits_{T^*_xM} \int\limits_{B_\delta(x)}
  \chi(x,y) e^{i \varphi_H(\alpha_x,x)} u(y) | \det \exp_* |
  dy d\xi ,
\ee
where $\chi$ is a bump function as in (5.26) and $\varphi_H$ is the 
global phase function of $\Gamma(f)$ defined in Lemmas 5.13, 5.14. 
The formula (5.59)
holds modulo ${\cal U}\Psi^{-\infty}(M)$.

We define the local section $\sigma$ as follows: 
Let $f \in {\cal U} \cup {\cal D}^\infty_{\theta,0}(T^*M \setminus 0)$ 
and define $\sigma (f) \in ({\cal U} F^{0,k}(\infty))_*$ by
\be\label{section}
  \sigma(f) u(x) := (2\pi)^{-n} \int\limits_{T^*_xM}
  \int\limits_{B_\delta (x)} \chi(x,y) e^{i \varphi_H (\alpha_x,y)}
  u(y) |\det \exp_*| dy d\xi ,
\ee
where $\varphi_H(\alpha_x,y)=\varphi_0(\alpha_x,y)+H(\alpha_x)$,
$H=\Psi(f)$, i. e. 
$\varphi_H(\alpha_x,y)=\langle\alpha_x,\exp^{-1}_x(y)\rangle-
\langle\alpha_x, \exp^{-1}_x(\pi(f^{-1}(\alpha_x)))\rangle$.

The operator $\sigma(f)$ is a FIO with smooth phase function
$\varphi_{\Psi(f)}$ and amplitude $a=1$. 
Moreover, $\sigma(f)$ is invertible
modulo smoothing operators  since $f$ is invertible and its principal symbol
is $a=1$, hence 
$\sigma(f) \in ({\cal U} F^{0,k}(\infty))_*$ for any $k$. 
Furthermore, $\pi \sigma (f)=f$, hence  $\sigma$ is a local
section of the exact sequence (5.24).

\setcounter{equation}{0}
\section {$({\cal U} \Psi^{0,k})_*$ as ILH Lie group}
                                               
We want to endow $({\cal U} \Psi^{0,k})_*$ with the structure of an ILH 
Lie group. Consider $[A] \in {\cal U} \Psi^q$ (which means
${\cal U} \Psi^q / {\cal U} \Psi^{-\infty}$) with principal symbol
$a_q$ which is globally defined and behaves well under transformations.
Let $A_q$ be the operator given by (5.26) with total symbol $a_q$.
Then $[A-A_q] \in {\cal U} \Psi^{q-1}$. Let $a_{q-1}$ be the 
principal symbol of $A-A_q$ and $A_{q-1}$ given by (5.26).
Then $[A-A_q-A_{q-1}] \in {\cal U} \Psi^{q-2}$. Continuing in the
manner, we obtain an assignment
\be   
  [A] \in {\cal U} \Psi^q \longmapsto 
  (a_q(x,\xi), a_{q-1}(x,\xi), \dots ) 
\ee   
where $a_{q-1} \in C^\infty (T^*M \setminus 0)$, satisfies (5.25) and
is homogeneous of degree $q-1$ in $\xi$, hence not square integrable.
We consider their restriction to $S(T^*M)$ and denote 
$a_{q-j}|_S$ again by $a_{q-j}$. The map (6.1) is still a vector space
isomorphism. Fix some $k$ and consider the assignment
\be    
  [A] \in {\cal U} \Psi^{q,k} \longmapsto
  (a_q, a_{q-1}, \dots , a_{-k}) 
\ee    
We introduce a uniform Sobolev topology on $ {\cal U} \Psi^{q,k}$.
Let $\delta>0, s>n$ and set
\bea   
V_\delta & = & \Big\{ ([A], [A']) \in {\cal U} \Psi^{q,k}\times {\cal U}
\Psi^{q,k} ~\Big|~|[A]-[A']|^2_{q,k,s} := \nonumber \\
   &: = & |a_q-a'_q|^2_{q+k+s} + |a_{q-1}-a'_{q-1}|^2_{q+k+s-1} 
   + \dots + |a_{-k}-a'_{-k}|^2_{s} < \delta^2 \Big\}   \nonumber
\eea   
       
\me    
\no    
{\bf Lemma 6.1} {\it 
$L=\{V_\delta\}$ is a basis for a metrizable uniform structure
${\cal A}^{q,k,s}$.
}      
       
\me    
\no    
We omit the very simple proof.
\hfill $\Box$
                               
Let ${\cal U} \Psi^{q,k,s} = 
\overline{{\cal U} \Psi^{q,k}}^{||_{q,k,s}}$
be the completion. 
       
\me    
\no    
{\bf Proposition 6.2.} {\it
${\cal U} \Psi^{q,k,s}$ is the topological sum of its arc
components,
\be    
  {\cal U} \Psi^{q,k,s} = \sum\limits_{i \in I} comp ([A_i]) ,
\ee    
and each component is a smooth Hilbert manifold. Here
\be     
  comp([A])= \Big\{ [A'] \in {\cal U} \Psi^{q,k,s} \Big| 
  |[A]-[A']_{q,k,s}| < \infty \Big\} .
\ee     
}       
        
\me     
\no     
We omit the simple proof which is performed for spaces of connections
or spaces of metrics e. g. in \cite{Ei1}, \cite{Ei4}. 
\hfill $\Box$
        
\no      
In ${\cal U}\Psi^{0,k}$ composition is well defined.
            
\me         
\no        
{\bf Proposition 6.3.} {\it
Composition in ${\cal U} \Psi^{0,k}$ extends continuously to 
${\cal U} \Psi^{0,k,s}$. More precisely, composition is a continuous 
map                    
\be                   
  comp([A]) \times comp([B]) = comp([A \circ B]) .
\ee
}

\me
\no                                                  
{\bf Proof.} The elements of ${\cal U} \Psi^{0,k}$ are dense in  
${\cal U} \Psi^{0,k,s}$. Fix $[A], [B], [A \circ B] \in  
{\cal U} \Psi^{0,k}$. We have to show that for 
$[A'] \in comp([A]), [B'] \in comp([B])$
\be
  [A' \circ B'] \in comp ([A \circ B])
\ee
and that this map is continuous.
The proof of (6.6) will also include the proof of continuity.
Represent the operators $[A], [B], [A \circ B]=[C], [A'], [B'], 
[A' \circ B']=[C']$ by the symbols

$a=a(x,\xi)=a_0+a_{-1}+ \dots + a_{-k}$, 

$b=b(x,\xi)=b_0+b_{-1}+ \dots + b_{-k}$, 

$c=c(x,\xi)=c_0+c_{-1}+ \dots + c_{-k}$, 

$a'=a'(x,\xi)=a'_0+a'_{-1}+ \dots + a'_{-k}$, etc..

\no
Then
\[c_0=a_0 b_0, \en c'_0=a'_0 b'_0 \]
\bea
    c_{-1} &=& a_0b_{-1}+a_{-1}b_0 + \sum\limits^n_{i=1} 
    \partial_{x_i} a_0 \partial_{\xi_i} b_0 , \nonumber \\ 
    c'_{-1} &=& a'_0b'_{-1}+a'_{-1}b'_0 + \sum\limits^n_{i=1} 
    \partial_{x_i} a'_0 \partial_{\xi_i} b'_0 , \nonumber \en \mbox{etc.}
\eea    
Then
\bea
    |a_0b_0 - a'_0b'_0|_{k+s} &=&  
    |a_0b_0 - a'_0b_0 + a'_0b_0 - a'_0b'_0|_{k+s} \nonumber \\
    &\le& |(a_0-a'_0)b_0|_{k+s} + |a'_0(b_0-b'_0)|_{k+s} \nonumber \\ 
    &\le& |(a_0-a'_0)b_0|_{k+s} + |(a'_0-a_0)(b_0-b'_0)|_{k+s}
    + |a_0(b_0-b'_0)|_{k+s} \nonumber \\
    &\le& |a_0-a'_0|_{k+s} {}^{b,k+s}|b_0| + C_0 |a'_0-a_0|_{k+s}
    |b_0-b'_0|_{k+s} + {}^{b,k+s}|a_0| \cdot |b_0-b'_0|_{k+s} 
    . \nonumber \\
    {}
\eea 
Here we applied  the module structure theorem in the middle term and in the
boundary terms (5.25) for $a_0, b_0$. Next we have to estimate
\be
   |c_{-1}-c'_{-1}|_{k-1+s} \le |a_0b_{-1}-a'_0b'_{-1}|_{k-1+s}
   + |a_{-1}b_0-a'_{-1}b'_0|_{k-1+s} + \sum\limits^n_{i=1}
   |\partial_{x_i} a_0 \partial_{\xi_i} b_0 - 
   \partial_{x_i} a'_0 \partial_{\xi_i} b'_0|_{k-1+s}
\ee
Each single term in (6.8) can be estimated as in (6.7). The general case can 
be proved by a simple but very extensive induction. The key
inequality is (6.7) (applied with other indices).
           
\hfill $\Box$
          
\me       
\no       
{\bf Proposition 6.4.} {\it
$({\cal U} \Psi^{0,k,s})_*$ is a Hilbert Lie group.
}        
        
\me     
\no     
{\bf Proof.} According to (6.4), each component is an affine Hilbert
space. Proposition 6.3 implies that ${\cal U} \Psi^{0,k,s}$ is an affine
Hilbert algebra. Then $({\cal U} \Psi^{0,k,s})_*$  is open and
$A \mapsto A^{-1}$ is a homoemorphism. Composition in 
${\cal U} \Psi^{0,k,s}$ is even smooth since it is bilinear and 
continuous. Hence the same holds in $({\cal U} \Psi^{0,k,s})_*~.$
\hfill $ \Box$   
        
\me   
\no      
{\bf Remarks.} 
      
\no      
{\bf 1.} If the underlying manifold $M$ is compact then $({\cal U} 
\Psi^{0,k,s})_*$ consists of one component (as 
${\cal U} \Psi^{0,k,s}$ does). On open manifolds, 
$({\cal U} \Psi^{0,k,s})_*$ consists of uncountably many components
(as ${\cal U} \Psi^{0,k,s}$ does). 
          
\no   
{\bf 2.} If we look in $({\cal U} \Psi^{0,k,S})_*$ for those $As$
which are invertible in their own component then this component must be
$comp([I])$. This can be proved by easy calculations and estimates .

\hfill $\Box$
           
\me        
We established the following inverse systems                      
\[ \cdots \longrightarrow {\cal U} \Psi^{0,k,s+1}
   \longrightarrow {\cal U} \Psi^{0,k,s} \longrightarrow \cdots
   \longrightarrow {\cal U} \Psi^{0,k,s_0} \] 
and                   
\[ \cdots \longrightarrow ({\cal U} \Psi^{0,k,s+1})_* 
   \longrightarrow ({\cal U} \Psi^{0,k,s})_* \longrightarrow \cdots
   \longrightarrow ({\cal U} \Psi^{0,k,s_0})_* \] 
with                              
${\cal U} \Psi^{0,k} = \lim\limits_\leftarrow {\cal U} \Psi^{0,k,s}$ and 
$({\cal U} \Psi^{0,k})_* = \lim\limits_\leftarrow 
({\cal U} \Psi^{0,k,s})_*$. \\    
The tangent space              
$T_{[I]} ({\cal U} \Psi^{0,k,s})_*$ can be described as follows:
\bea        
   T_{[I]} ({\cal U} \Psi^{0,k,s})_* & = & 
   T_{[I]} (comp([I]))_* = T_{[I]} comp([I]) \nonumber \\ 
   & = & T_{[I]}([I]+ \{ [A] \in {\cal U} \Psi^{0,k,s} | 
   |a_0|^2_{k+s} + |a_{-1}|^2_{k-1+s} + \dots + 
   |a_{-k}|^2_s < \infty \} ) \nonumber \\
   & = & \{ [A] \in {\cal U} \Psi^{0,k,s} | 
   |a_0|^2_{k+s} + |a_{-1}|^2_{k-1+s} + \dots + 
   |a_{-k}|^2_s < \infty \} \nonumber \\
   & = & \mbox{zero component of} \en {\cal U} \Psi^{0,k,s} \nonumber \\
   & \cong & \Omega^{0,2,k+s} (S) \oplus \dots \oplus
   \Omega^{0,2,s} (S) . \nonumber
\eea         
This yields an inverse system
\[ \cdots \longrightarrow T_{[I]} ({\cal U} \Psi^{0,k,s+1})_*  
   \longrightarrow T_{[I]} ({\cal U} \Psi^{0,k,s})_* 
   \longrightarrow \cdots  \]                                         
and             
\[ \lim\limits_{\leftarrow} T_{[I]} ({\cal U} \Psi^{0,k,s})_* =
   \{ [A] \in {\cal U} \Psi^{0,k} | a_0, a_{-1}, \dots a_{-k} 
   \in \Omega^{0,2,s} (S) \en \mbox{for all} \en s \} . \]
We denote the latter space by $L {\cal U} \Psi^{0,k} :=
\lim\limits_{\leftarrow} T_{[I]} ({\cal U} \Psi^{0,k,s})_*$.   
This is in fact a Lie algebra with respect to the
bracket of $\Psi$DOs. 
                   
\me                   
Hence we proved the following 
                       
\me            
\no                     
{\bf Theorem 6.5.} {\it
Assume $(M^n,g)$ open satisfying the conditions $(I)$ and $(B_\infty)$.
Then     
$\{ ({\cal U} \Psi^{0,k})_* , ({\cal U} \Psi^{0,k,s})_* | s>n \}$
is an ILH Lie group and each $({\cal U} \Psi^{0,k,s})_*$ is a
smooth Hilbert Lie group. Its ILH Lie algebra is 
$\{ L {\cal U} \Psi^{0,k}, T_{[I]} ({\cal U} \Psi^{0,k,s})_* | s>n \}$.
Here $L {\cal U} \Psi^{0,k}$ is isomorphic, as topological vector space, to
$\sum\limits^k_1 \Omega^{0,2,\infty} (S)$, where 
$\Omega^{0,2,\infty} (S) = \bigcap\limits_{S}
\Omega^{0,2,S} (S)$.    
}                           
\hfill $\Box$           
                 
\me              
\no              
{\bf Corollary 6.6.} {\it
$({\cal U} \Psi^0)_*$ has the structure of a direct limit of ILH Lie
groups,           
\[ ({\cal U} \Psi^0)_* = \lim\limits_{\stackrel{\longrightarrow}{k}}~
   ({\cal U} \Psi^{0,k})_* ~. \]
}                 
\hfill $\Box$

\setcounter{equation}{0}
                   
\section{An ILH Lie group structure for invertible Fourier integral
operators}      

We consider our exact sequence (5.24), 
\[ I \longrightarrow ({\cal U} \Psi^{0,k})_* \longrightarrow 
   ({\cal U} F^{0,k} (\infty))_* \longrightarrow 
   {\cal D}^\infty_{\theta,0} \longrightarrow e \]
and perform the 4th steps of our program described  after Lemma 5.9. We know
already that $({\cal U}\Psi^{0,k})_*$ and ${\cal D}^\infty_{\theta,0}$ are
ILH Lie groups, in particular they are topological groups.

We first consider an exact sequence of abstract groups
\be 
  I \longrightarrow {\cal H} 
  \begin{array}{c} j \\[-2ex] \longrightarrow \\[-2ex] {} \end{array} 
  {\cal G}                                                    
  \begin{array}{c} \pi \\[-2ex] \longrightarrow \\[-2ex] {} \end{array} 
  {\cal Q} \longrightarrow e
\ee                
We assume that ${\cal H}$ and ${\cal Q}$ are topological groups and
construct a topology on ${\cal G}$ such that (7.1) becomes an exact sequence
of topological groups. 
In a second step we will sharpen the construction to
the case of ILH Lie groups.
The frame work of this approach is given in Adams-Ratiu-Schmid \cite{ARS2}. We
recall without proofs the  facts established there and concentrate our
attention to the new  features coming from the openness of the underlying
manifold $M$.

\me
\no
{\bf Lemma 7.1.} {\it 
Let
\be
  I \longrightarrow H 
  \begin{array}{c} j \\[-2ex] \longrightarrow \\[-2ex] {} \end{array} 
  G                                                     
  \begin{array}{c} \pi \\[-2ex] \longrightarrow \\[-2ex] {} \end{array} 
  Q \longrightarrow e
\ee   
be an exact sequence of groups where $H$ and $Q$ are topological groups.
Let $U \subset Q$ be a neighborhood of the identity $e\in Q$ and  
$\sigma: U \longrightarrow G$ a local section of $\pi$; let $V \subset U$
be a neighborhood of $e \in Q$ such that $V \cdot V^{-1} \subset U$
and assume               
                           
\no                       
A) the map $V \times V \times H \longrightarrow H$  given by
\[ (f_1,f_2,h) \longrightarrow \sigma (f_1) \sigma (f_2)^{-1}
   h \sigma (f_1 f^{-1}_2)^{-1} \]
is continuous;             
                           
\no                        
B) for each $g \in G$ and $W \subset U$ such that 
$\pi(g) W \pi(g)^{-1} \subset U$ , the map $W \times H \longrightarrow H$
given by                    
\[ (f,h) \longrightarrow g h \sigma(f) g^{-1} \sigma(\pi (g) f 
    \pi(g)^{-1})^{-1} \]    
is continuous.               
                             
\no                          
Then $G$ can be made into a topological group such that $j, \pi$ and $\sigma$
are continuous and $\pi$ is open.
}                            
\hfill $\Box$                
                             
\me                           
\no                           
{\bf Remark.} {\it             
If $Q$ is connected then the condition $A)$ is already sufficient. }
\hfill $\Box$                  

\me
Assume now ${\cal H}$ and ${\cal Q}$ in (7.1) to be ILH Lie groups  with 
${\cal H} = \lim\limits_{\leftarrow} H^s$  and  
${\cal Q} = \lim\limits_{\leftarrow} Q^t$, $s \ge s_0$, 
$t \ge t_0$. 
Denote by ${\cal H}^s$ and ${\cal Q}^t$ the spaces
${\cal H}$ and ${\cal Q}$ with the coarser topologies of 
${\cal H}^s$ and ${\cal Q}^t$, respectively.
Then ${\cal H}^s$ and ${\cal Q}^t$ are topological groups.

Suppose $ U \subset {\cal Q}^{t_0}$ open (hence open in all 
${\cal Q}^{t}$) , $ \sigma: U \longrightarrow {\cal G}$ a
local section and $V \subset U$ open in 
${\cal Q}^{t_0}$ such that $V \cdot V^{-1} \subset U$.
Denote $V$ with the topology from ${\cal Q}^t$ by $V^t$. Assume that 
for each $t \ge t_0$ there exists an $s(t)$ such that

\no
A') The map $V^t \times V^t \times {\cal H}^{s(t)} \longrightarrow 
{\cal H}^{s(t)}$ defined by  
\[ (f_1,f_2,h) \longrightarrow \sigma(f_1) \cdot \sigma(f_2)^{-1}
   h \sigma (f_1f^{-1}_2)^{-1} \]
is continuous;

\no
B') for each $g \in {\cal G}$ and $W \subset U$ with 
$\pi(g) W \pi(g)^{-1} \subset U$, the map 
$W^t \times {\cal H}^{s(t)} \longrightarrow 
{\cal H}^{s(t)}$ defined by  
\[ (f,h) \longrightarrow g h \sigma(f) g^{-1} \sigma(\pi(g) f 
   \pi(g)^{-1})^{-1} \]
is continuous.

Then we derive from Lemma 7.1. that for each $t$ the group ${\cal G}$ becomes
a topological group ${\cal G}^t$, and we obtain an exact sequence of
topological groups for all $t \ge t_0$, 
\be 
  I \longrightarrow {\cal H}^{s(t)} 
  \begin{array}{c} j \\[-2ex] \longrightarrow \\[-2ex] {} \end{array} 
  {\cal G}^t
  \begin{array}{c} \pi \\[-2ex] \longrightarrow \\[-2ex] {} \end{array} 
  {\cal Q}^t \longrightarrow e
\ee                
Taking the inverse limit topology, we obtain (7.1) as a sequence of 
topological groups.

We recall the right uniform structure of a topological group $G$.
A sequence 
$(x_n)_n$ is a Cauchy sequence with respect to this structure if for any
neighborhood $V$ of the identity there exists an $n_0$ s. t. 
$x_n x^{-1}_m \in V$ for all $m,n \ge n_0$.

\me  
\no
{\bf Lemma 7.2.} {\it
Let $X$ be a locally Hilbert topological group. Then $X$ is complete in
its right uniform structure.
}
\hfill $\Box$

\me
\no
{\bf Corollary 7.3.} {\it
$\overline{{\cal H}^s} = H^s$ and 
$\overline{{\cal Q}^t} = Q^t$.
}

\me
\no
{\bf Proof.} 
${\cal H}^s$ is dense in $H^s$, ${\cal Q}^t$ is dense in $Q^t$. 
Then apply Lemma 7.2.
\hfill $\Box$

\me
Let $G^t$ be the completion of ${\cal G}^t$ with respect to the right
uniform structure. It is not yet clear that $G^t$ is  
a topological group or even a group. 

\me
\no
{\bf Lemma 7.4.} {\it 
Let $\tilde{U}^{t_0} \subset Q^{t_0}$ be open and 
$\tilde{U}^{t_0} \cap {\cal Q}^{t_0} = U^{t_0}$. For
$t \ge t_0$ , let $\tilde{U}^{t} = \tilde{U}^{t_0} \cap Q^t$ and
$\tilde{V}^{t_0} \subset \tilde{U}^{t_0}$ open with  
$\tilde{V}^{t_0} \cap {\cal Q}^{t_0} = V^{t_0}$ and 
$\tilde{V}^{t} = \tilde{V}^{t_0} \cap Q^t$ , $t \ge t_0$.
Assume that $\sigma: U \longrightarrow {\cal G}$
extends to a local section 
$\tilde{\sigma} : \tilde{U}^t \longrightarrow G^t$ and that
            
\no          
A'') The map $V^t \times V^t \times {\cal H}^{s(t)} \longrightarrow
{\cal H}^{s(t)}$ from A') actually extends by $\tilde{\sigma}$
to a continuous map 
$\tilde{V}^t \times \tilde{V}^t \times H^{s(t)} \longrightarrow
H^{s(t)}$.  \\   
Then $G^t$ is a topological group.
}                
\hfill $\Box$    
                  
Next we want to endow $G^t$ with a manifold structure. Define
$\Psi_I : \pi^{-1}(U^t) \longrightarrow U^t \times H^{s(t)}$ by
\be                
  \Psi_I : g \longmapsto (\pi(g), g  \sigma(\pi(g))^{-1}) .
\ee
We try to consider this as a bundle chart and move this around on 
$G^t$ by right translations: Let $g_0 \in G^t$, $q_0=\pi(g_0)$
and define on $\pi^{-1}(U q_0)$
\[ \Psi_{g_0} : \pi^{-1} (U^t q_0) \longrightarrow 
   U^t q_0 \times H^{s(t)}  \]
by                   
\be                  
 \Psi_{g_0} : g \longmapsto (\pi(g), g g^{-1}_0 \sigma ( \pi(g)
q^{-1}_0)^{-1} )
\ee                  
To obtain an atlas, we need the transition condition that
\be                  
  \Psi_{\tilde{g}_0} \circ \Psi^{-1}_{g_0} :
  (U^t q_0 \cap U^t \tilde{q}_0) \times H^{s(t)} 
  \longrightarrow (U^t q_0 \cap U^t \tilde{q}_0)
  \times H^{s(t)} 
\ee                   
is $C^{k(t)}$. According to our definitions this is the map
\be                    
  (f,h) \longmapsto (f, h \sigma (f q^{-1}_0) g_0 \tilde{g}^{-1}_0
  \sigma ( f \tilde{q}^{-1}_0)^{-1})  .
\ee                    
Hence we add the condition
                       
\no                    
C) For $g_0, \tilde{g}_0 \in G^t$ with $\pi(g_0) = q_0$, 
$\pi(\tilde{g}_0) = \tilde{q}_0$ assume that the map (7.6), (7.7) is
$C^{k(t)}$, where $k(t)$ is an increasing function of $t$.
                        
Finally to construct an ILH Lie group structure for ${\cal G}$ we
need that multiplication $G^{r+k} \times G^r \longrightarrow G^r$
is $C^k$ for $k \le k(r)$. As pointed out in \cite{ARS2}, p. 30, this leads 
to the final condition  
                        
\no                     
D) Assume for $a \in G^{r+k}$, $b \in G^r$ with $\alpha=\pi(a)$,
$\beta=\pi(b)$ that the map            
\[ U^{r+k} \alpha \times U^r \beta \times 
   H^{s(r)} \longrightarrow H^{s(r)} , \]
\[ (f_1,f_2,h) \longmapsto \sigma(f_1 \alpha^{-1} a h 
   \sigma(f_2 \beta^{-1}) a^{-1} \sigma(f_1, f_2 \beta^{-1} 
   \alpha^{-1})^{-1}) \]
is $C^k$ as long as $k \le k(r)$.
                         
\me                      
Summarizing, we obtained 
                         
\setcounter{equation}{0} 
\me                      
\no                      
{\bf Theorem 7.5} {\it Let  
\be                      
  I \longrightarrow {\cal H} \longrightarrow {\cal G}
  \longrightarrow {\cal Q} \longrightarrow e
\ee                      
be an exact sequence of groups where ${\cal H}$ and ${\cal Q}$ have
ILH Lie group structures. Suppose there exists a local section 
$\sigma:  U \longrightarrow {\cal G}$ satisfying the conditions
A''), B'), C) and D). Then ${\cal G}$ has an ILH Lie group structure
and (7.1) becomes an exact sequence of ILH Lie groups. If ${\cal Q}$
is connected then the condition B') follows from A'').      
}                                       
\hfill $\Box$  

For the proofs of 7.1 - 7.5 we refer to
Adams-Ratiu-Schmid \cite{ARS2}.                   
                                
\me                            
Now we apply Theorem 7.5 to our situation where 
${\cal H} =\lim\limits_{\leftarrow} \{H^s | s \ge s_0 \}$ with 
$H^s = ({\cal U } \Psi^{0,k,s})_*$ , $s_0=n+1$ ,  
$~{\cal Q} =\lim\limits_{\leftarrow} \{Q^t | t \ge t_0 \}$ with  
$Q^t={\cal D}^t_{\theta,0}$, $t_0=n+1$ and 
${\cal G} = ({\cal U} F^{0,k}(\infty))_*$,
\[ I \longrightarrow ({\cal U } \Psi^{0,k})_* \longrightarrow 
   ({\cal U} F^{0,k}(\infty))_* \longrightarrow 
   {\cal D}^\infty_{\theta,0} (T^*M \setminus 0) 
   \longrightarrow e\]
\[ \sigma : {\cal U} \subset {\cal D}^\infty_{\theta,0}
   (T^*M \setminus 0) \longrightarrow ({\cal U} F^{0,k}(\infty))_*  \]
defined by (\ref{section}). We have to verify the conditions above.
         
As always now, we assume $(M^n,g)$ with $(I)$, $(B_\infty)$ and
$\inf \sigma_e (\bigtriangleup_1(S) 
|_{ker \bigtriangleup_1(S))^\perp } ) > 0$.
             
\setcounter{equation}{7}

\me            
\no             
{\bf Proposition 7.6.} {\it
Let $\sigma : {\cal U} \subset {\cal D}^\infty_{\theta,0} 
(T^*M \setminus 0) \longrightarrow ({\cal U} F{0,k}(\infty))_*$
be defined by (\ref{section}), $V \subset {\cal U}$ be a neighborhood of
$e \in {\cal D}^\infty_{\theta,0}$ such that $V \cdot V^{-1} \in {\cal U}$.
Then the condition A) is satisfied, i. e. the map
\be                    
  V \times V \times ({\cal U} \Psi^{0,k})_* \longrightarrow
  ({\cal U} \Psi^{0,k})_* ,
\ee 
\be                       
  (f_1,f_2,A) \longrightarrow \sigma (f_1) \sigma(f_2)^{-1} A 
  \sigma (f_1 f^{-1}_2)^{-1} ,
\ee
is continuous and extends as $C^r$ map to certain Sobolev completions,
which will be specified below.
}                             
                              
\me                           
\no                           
{\bf Proof.}                   
We are done if we can prove the following fact. Let 
$A, B \in ({\cal U} F^{0,k})_*$ near the identity. Then the 
symbols of $A \circ B$ and $A^{-1}$ depend continuously on the 
symbols and phase functions of $A$ and $B$, i. e. if 
$\varphi_{H_1}, H_1=\Psi(f_1)$, $\varphi_{H_2}, H_2=\Psi(f_2)$,    
$\varphi_H, H=\Psi(f_1f_2)$, $\varphi_{H^-}, H^-=\Psi(f^-_1)$,                                                                     
are global phase functions for $A, B, A \circ B, A^{-1}$, 
respectively, and $a, b, c, c^-$ are global symbols, 
$A', B'$ operators of the same kind with
$\varphi_{H'_1}, \varphi_{H'_2}, \varphi_{H'} \varphi_{H'^-}, a'$
etc., then                     
\bea                           
  |H-H'|_{Sob} & \le & P_1 (|H_1 - H'_1|_{Sob}, 
  |H_2 - H'_2|_{Sob}) , \\  
  |H^- - H'^-|_{Sob} & \le & P_2 (|H_1 - H'_1|_{Sob}) , \\
  |c-c'|_{Sob} & \le & P_3 (|a-a'|_{Sob}, |b-b'|_{Sob}, 
  |H_1 - H'_1|_{Sob}, |H_2 - H'_2|_{Sob}) , \\  
  |a^- - a'^-|_{Sob} & \le & P_4 (|a-a'|_{Sob}, 
  |H_1 - H'_1|_{Sob}) .   
\eea                             
Here the $P_j$ are polynomials without constant terms and 
$|\en|_{Sob}$ means certain Sobolev norms, the Sobolev index of which
remains still open for a moment. Cover $M$ by a uniformly locally finite
cover $U = \{ U_i \}_i$ of normal charts. Then it is a
well known fact that there exist constants $C_1, C_2$
s. t.                             
\be                               
  C_1 \sum\limits_i |\en|^2_{Sob,U_i} \le
  |\en|^2_{Sob} \le         
  C_2 \sum\limits_i |\en|^2_{Sob,U_i} .                     
\ee                               
(7.14) immediately implies that we are done if we can show (7.10)--(7.13)
locally, i. e.                    
\bea                              
  |H-H'|_{Sob,  U_i} & \le & C \cdot P_1 
  (|H_1 - H'_1|_{Sob, U_i}, 
  |H_2 - H'_2|_{Sob, U_i}) , \\  
  |H^- - H'^-|_{Sob, U_i} & \le & 
  C \cdot P_2 (|H_1 - H'_1|_{Sob, U_i}) , \\
  |c-c'|_{Sob, U_i} & \le & C \cdot P_3 ( \dots ) , \\ 
  |a^- - a'^-|_{Sob,U_i} & \le & C \cdot P_4 ( \dots ). 
\eea                                                                           
with $C$ independent of $i$. This is more or less explicitly done in 
\cite{ARS2}, lemma 4.2 and 4.3, p. 32--35. We recall the initial step. In
local coordinates we can write   
\bea
   A u(x) & = & (2 \pi)^{-n} \int \int e^{i(x-y)\xi+H_1(x,\xi)}
   a(x,\xi) u(y) dy d\xi, \nonumber \\                                          
   B u(x) & = & (2 \pi)^{-n} \int \int e^{i(x-y)\xi+H_2(x,\xi)}
   b(x,\xi) u(y) dy d\xi, \nonumber \\   
   A \circ B u(x) &  = & (2 \pi)^{-n} \int \int 
   e^{i(x-y)\eta+H_1(x,\eta)+y\xi+H_2(y,\xi)}
   a(x,\eta) b(y,\eta) \hat{u}(\xi) d\xi dy d\eta = \nonumber \\
   & = & (2 \pi)^{-n} \int \int e^{i(x\xi+H(x,\xi))}
   c(x,\xi) \hat{u}(\xi) d\xi \nonumber                                          
\eea   
with $H_1=\Psi(f_1), H_2=\Psi(f_2), H=\Psi(f_1 \circ f_2)$.
Hence  
\be                       
  c(x,\xi) \sim (2 \pi)^{-n} \int \int e^{i\varphi(x,y,\xi,\eta)}
  a(x,\eta) b(y,\xi) dy d\eta 
\ee                       
with                      
$\varphi(x,y,\xi,\eta)=(x-y) \cdot (\eta-\xi) + H_1(x,\eta) +
H_2(y,\xi) - H(x,\xi)$.   
Locally we have for (7.17) to establish the continuity of the map
\be                       
  S^{0,k} \times S^{0,k} \times {\cal W} \times {\cal W} 
  \longrightarrow S^{0,k} 
\ee                          
\be                       
  (a(x,\xi), b(x,\xi), H_1, H_2) \longmapsto c(x,\xi)
\ee                       
and its extension to a certain Sobolev completion, 
\be                       
  S^{0,k,\tilde{s}(t)} \times S^{0,k,\tilde{s}(t)}
  \times {\cal W}^{t+r} \times {\cal W}^{t+r} \longrightarrow
  S^{0,k,\tilde{s}(t)} .  
\ee    
Here $S^{o,k}$ is defined by ${\cal U} S^0 / {\cal U} S^{-k-1}$,
${\cal U} S^q = {\cal U} S^q (B) \equiv 
{\cal U} S^q (B \times \R^n)$. 
If we could establish the continuity of (7.20)--(7.22), 
then we would have (7.17) by difference constructions with the
gotten formulas, if these formulas permit such a construction. This is
in fact the case. We refer to \cite{ARS2}, p. 34--35. The main point is to
calculate or estimate (7.19). This has been done in \cite{ARS2} by the method
of stationary phase (as one would expect). Finally the Sobolev index
in (7.17) and (7.20) is $\tilde{s}(t)=2(t-k-1)$. The estimate (7.15) is
trivial as we have see from the last representation for $A \circ B$. The same
holds for (7.16) and (7.18) which follows from Lemma 4.3 in \cite{ARS2} after
some calculations.
\hfill $\Box$            
         
\me                      
Now Lemma 7.1 and Proposition 7.6 imply the following
                       
\me      
\no                   
{\bf Theorem 7.7.} {\it
$({\cal U} F^{0,k} (\infty) )_*$ is a topological group , (5.24) is an
exact sequence of topological groups and the local section $\sigma$
is continuous. }  
\hfill $\Box$     
          
\me               
We write in the sequel ${\cal U} F^{0,k}$ instead of
${\cal U} F^{0,k} (\infty)$ since we consider only that space.
Now let $ H^s = ({\cal U} \Psi^{0,k})^s_*$ and 
$Q^t = ({\cal D}^\infty_{\theta,0} (T^*M \setminus 0) )^t$
the space ${\cal D}^\infty_{\theta,0}$ with coarser topology of
${\cal D}^t_{\theta,0} (T^*M \setminus 0)$. 
Proposition 7.6 implies that for $t>2n$ and $s(t)=t-2(k+1)>n$ the map
\[ V^t \times V^t \times ({\cal U} \Psi^{0,k})^{s(t)}_*  
   \longrightarrow  ({\cal U} \Psi^{0,k})^{s(t)}_*  \]
\[ (f_1,f_2,A) \longmapsto \sigma(f_1) \sigma(f_2)^{-1} 
   A \sigma (f_1 f^{-1}_2)^{-1} , \]         
is at least continuous. Hence we obtain from Proposition 7.1 that 
$({\cal U} F^{o,k})_*$ becomes a topological group 
$({\cal U} F^{o,k})^t_*$ s. t. 
\[ I \longrightarrow ({\cal U} \Psi^{0,k})^{s(t)}_*
   \begin{array}{c} j \\[-2ex] \longrightarrow \\[-2ex] {} \end{array}
   ({\cal U} F^{0,k})^t_* 
   \begin{array}{c} \pi \\[-2ex] \longrightarrow \\[-2ex] {} \end{array}
   ({\cal D}^\infty_{\theta,0} (T^*M \setminus 0) )^t
   \longrightarrow e \]   
is an exact sequence of topological groups for 
$t \ge t_0= \max \{2n, n+2(k+1)\}$. Let $({\cal U} F^{0,k,t})_*$
be the completion with respect to the right uniform structure. We will
show that this is a topological group. For this we have to show that
the local section $\sigma$ extends to a local section 
\[ \tilde{\sigma} : \tilde{{\cal U}}^t \longrightarrow 
   ({\cal U} F^{0,k,t})_* \]
and the map in condition A') extends to a continuous map A")
                        
\no                     
\[ \mbox{A'')} \quad \quad \quad \tilde{V}^t \times \tilde{V}^t \times 
({\cal U} \Psi^{0,k,S(t)})_*  \longrightarrow 
({\cal U} \Psi^{0,k,S(t)})_* .\] 
Consider first the extension $\sigma \longrightarrow \tilde{\sigma}$.
Let $f \in \tilde{{\cal U}}^t$. The $\tilde{\sigma}(f)$ should be an
FIO of order 0 with  symbol 
$a(x,\xi)=\sum\limits^k_{j=0} a_{-j}(x,\xi)$ 
where                   
$(a_{-j} - a^\infty_{-j})|_{S} \in \Omega^{0,2,k+s(t)}(S)$
for a smooth $a^\infty_{-j} \in {\cal U} S^{-j}$ and phase
function $\varphi_h$ generated by $H=\Psi(f)$ with  
$f \in {\cal D}^t_{\theta,0}$ and 
$H|_S \in \Omega^{0,2,k+1}(S)$.
The definition (5.26)  still makes sense as oscillatory integral of
$a(x,\xi)$. The
$H(x,\xi)$ can be differentiated enough times and $t>2n$ will be
sufficient for this. The continuity of the extension follows from
the procedures in Proposition 7.6, i. e. we have the condition A'') for 
$t>2n$ and $s(t)>n$. Hence we have established 
                   
\me               
\no
{\bf Theorem 7.8.} {\it
Assume $t>2n$ and $s(t)=t-2(k+1)>n$. Then 
\[ I \longrightarrow ({\cal U} \Psi^{0,k,s(t)})_*
   \begin{array}{c} j \\[-2ex] \longrightarrow \\[-2ex] {} \end{array}
   ({\cal U} F^{0,k,t})_* 
   \begin{array}{c} \pi \\[-2ex] \longrightarrow \\[-2ex] {} \end{array}
   ({\cal D}^t_\theta (T^*M \setminus 0) )^t
   \longrightarrow e  \]
is an exact sequence of topological groups such that $j$ and $\pi$ are 
continuous and $\pi$ is open. Moreover, (\ref{section}) defines a continuous
local section
\[ \tilde{\sigma}: {\cal U}^t \subset 
   {\cal D}^t_{\theta,0} (T^*M \setminus 0) \longrightarrow 
   ({\cal U} F^{0,k,t})_* . \]
}
\hfill $\Box$

\me
It remains to assure the conditions C) and D). 

\me
\no
{\bf Lemma 7.9.} {\it
Assume $t>2n$, $s(t)=t-2(k+1)>n$.
For $A_0, \tilde{A}_0 \in ({\cal U} F^{0,k,t})_*$ let
$\pi(A_0)=f_0$, $\pi(\tilde{A}_0)=\tilde{f}_0$. 
Then the condition C) is satisfied, i.e.the map 

\no
  
\[ {\cal U}^t \cdot f_0 \cap {\cal U}^t \cdot \tilde{f}_0 \times
   ({\cal U} \Psi^{0,k,s(t)})_* \longrightarrow
   ({\cal U} \Psi^{0,k,s(t)})_* , \]
\[ (f, A) \longmapsto A  \sigma (f f^{-1}_0) A_0 \tilde{A}^{-1}_0
   \sigma (f \tilde{f}^{-1}_0)^{-1} \]
is of class $C^t$.
}

\me
\no
{\bf Proof.}
According to Theorem 6.5, multiplication in $({\cal U} \Psi^{0,k,s})_*$
is smooth. Hence it suffices to show that 
\[ f \longmapsto \sigma (f f^{-1}_0) A_0 \tilde{A}^{-1}_0 
   \sigma (f \tilde{f}^{-1}_0)^{-1} \]
is $C^t$. For this it is sufficient to show that the map
\[ f \longmapsto \sigma (f) A \sigma (f f_A)^{-1} \]
is $C^t$ for $f$ and $f f_A$ near $id$, where we have set
$A= A_0 \tilde{A}^{-1}_0$, $f_A = \pi(A)$. But this follows 
from Proposition 7.6.
\hfill $\Box$

\me
Hence we have

\me
\no
{\bf Theorem 7.10.} {\it
Assume $t>2n$, $s(t)=t-2(k+1)>n$. Then 
$({\cal U} F^{0,k,t})_*$ is a Hilbert manifold of class $C^t$
modeled by 
$\Omega^{0,2,t+1}(S) \times ({\cal U} \Psi^{0,k,s(t)})_* $
}                                                              
\hfill $\Box$

\me
\no
{\bf Remarks.}

\no
1. Charts in $({\cal U} F^{0,k,t})_*$ are defined by right translation
of a chart at $I$. This implies automatically that right translation
in $({\cal U} F^{0,k,t})_*$ is $C^t$   . 

\no
2. As we have seen that $({\cal U} \Psi^{0,k,s(t)})_*$ consists of 
uncountable many components. The same  holds then also for
$j ({\cal U} \Psi^{0,k,s(t)})_*$.
\hfill $\Box$

\me
The last condition we have to verify is D).

\me
\no
{\bf Lemma 7.11.} {\it
Assume $t>2n$, $s(t)=t-2(k+1)>n$. Then for
$A \in ({\cal U} F^{0,k,t+r})_*$, $B \in ({\cal U} F^{0,k,t})_*$,
$f_A=\pi(A)$ $f_B=\pi(B)$.
Then the condition D) is satisfied, i.e. the map

\no

\[ {\cal U}^{t+r}_{f_A} \times {\cal U}^{t+r}_{f_B} \times 
   ({\cal U} \Psi^{0,k,s(t)})_* \longrightarrow 
   ({\cal U} \Psi^{0,k,s(t)})_* , \]
\[ (f_1,f_2,P) \longrightarrow \sigma (f_1 f^{-1}_A) A P \sigma
   (f_2,f^{-1}_B) A^{-1} \sigma (f_1 f_2 f^{-1}_B f^{-1}_A \]
is of class $C^l$ for $l=\min \{ r, t \}$.
}

\me
\no
{\bf Proof.}
We have that the map in condition A') (hence B') is $C^t$ and multiplication
in
${\cal D}^\infty_{\theta,0} (T^*M \setminus 0)$ is $C^r$ as a map
${\cal D}^{t+r}_{\theta,0} (T^*M \setminus 0) \times 
{\cal D}^t_{\theta,0} (T^*M \setminus 0) \longrightarrow 
{\cal D}^t_{\theta,0} (T^*M \setminus
0)$ .                                     
\hfill $\Box$                                   

\me 
Summarizing, we state our final main result.

\me                                    
\no                                    
{\bf Theorem 7.12.} {\it
Assume $(M^n,g)$ is an open Riemannian manifold satisfying the conditions
$(I)$ and $(B_\infty)$ of bounded geometry and the condition
\[ \inf \sigma_e (\bigtriangleup_1 (S(T^*M)),g_S 
   |_{ker \bigtriangleup_1)^\perp})>0, \]
Then for any $~ k \in {\bf Z}_+ $

\no
1. $\{ {\cal D}^\infty_{\theta,0} (T^*M \setminus 0), 
{\cal D}^r_{\theta,0} (T^*M \setminus 0) | r \ge n+1 \}$ 
is an ILH Lie group.

\no
2. $\{ ({\cal U} \Psi^{0,k})_* ~, ({\cal U} \Psi^{0,k,s})_*~ 
| s \ge n+1 \}$ is an ILH Lie group and each 
$({\cal U} \Psi^{0,k,s})_*$ is a smooth Hilbert Lie group.

\no
3. $\{ ({\cal U} F^{0,k})_*, ({\cal U} F^{0,k,t})_* 
| t > \max \{ 2n, n+2(k+1) \} \}$ is an ILH Lie group with the
following properties:

\no
a. $({\cal U} F^{0,k,t})_*$ is a $C^t$ Hilbert manifold modeled on
$\Omega^{0,2,t+1} (S (T^*M)) \times ({\cal U} 
\Psi^{0,k,t-2(k+1)})_*$.  
Each component of  $({\cal U} F^{0,k,t})_*$ is
modeled by
$\Omega^{0,2,t+1}(S) \oplus \Omega^{0,2,k+t-2(k+1)}(S) \oplus 
\dots \oplus \Omega^{0,2,t-2(k+1)}(S)$.

\no
b. The inclusion $({\cal U} F^{0,k,t+1})_* \hookrightarrow 
({\cal U} F^{0,k,t+1})_*$ is $C^t$. 

\no
c. The group multiplication 
$({\cal U} F^{0,k})_* \times ({\cal U} F^{0,k})_* \longrightarrow 
({\cal U} F^{0,k})_*$ 
extends to a $C^l$ map
\[ ({\cal U} F^{0,k,t+r})_* \times ({\cal U} F^{0,k,t})_* 
   \longrightarrow ({\cal U} F^{0,k,t})_* , \]   
$l=\min \{ r,t \}$.  

\no
d. The inversion 
$({\cal U} F^{0,k})_* \longrightarrow ({\cal U} F^{0,k})_*$ 
extends to a $C^l$ map
\[ ({\cal U} F^{0,k,t+r})_*  \longrightarrow ({\cal U} F^{0,k,t})_* , \] 
$l=\min \{ r,t \}$.  

\no
e. Right multiplication for $A \in ({\cal U} F^{0,k,t})_*$ is a 
$C^t$ map $R_A : ({\cal U} F^{0,k,t})_* \longrightarrow 
({\cal U} F^{0,k,t})_*$
}                                     
\hfill $\Box$                                

\me
We finish  with a few remarks concerning the corresponding Lie algebras
with differ slightly from the statements in \cite{ARS2}. First we present the
version of \cite{ARS2}, 5.2 for the open case.

\me
\no      
{\bf Lemma 7.13.} {\it
The Lie algebra of $({\cal U} F^{0,k,s})_*$ is $Comp_0 ~\widetilde{{\cal U}
\Psi}^{1,k,s}$ ,the 0--component of  the completed space $\widetilde{{\cal U}
\Psi}^{1,k,s}$ of formal  pseudodifferential operators of order one modulo
those of order 
$-k-1$ with pure imaginary principal symbol. The Lie bracket 
corresponds to the commutator bracket.}

\me
\no
{\bf Proof.} 
Let $c(t)$ be a $C^1$ curve in $({\cal U} F^{0,k,s})_*$ with
$c(0)=I$. We can write $c(t)=P(t) \sigma(\pi(c(t)))$, where $P(t)$ is
a $C^1$ curve in ${\cal U} \Psi^{1,k,s}$ such that $P(0)=I$ and
$\sigma(\pi(c(t)))$ has the local expression
\[ \sigma(\pi(c(t)))u(x)=(2\pi)^{-n} \int \int 
   e^{i(\varphi_0(x,y,\xi)+H_1(x,\xi))} u(y) dy d\xi \]   
with $\varphi(x,y,\xi)=\langle x-y,\xi \rangle$,
$H_0(x,\xi)=0$, $H_1(x,\xi) |_S \in \Omega^{0,2,k+s+1}(S)$.
Differentiation at $t=0$ yields $c'(0)=P'(0)+A$, 
where locally
\[ A u(x)=(2\pi)^{-n} \int \int 
   e^{i \varphi_0(x,y,\xi)} i H'_0(x,\xi) u(y) dy d\xi . \]
$H_t(x,\xi)$ is homogeneous of degree one in $\xi$ for each $t$. 
Hence the same holds for $H'_0(x,\xi)$, 
$H'_0 \in i \cdot {\cal H}^1 \Omega^{0,2,k+s+1} (T^*M \setminus 0)$.
This implies 
$A \in i \cdot Comp_0~ \widetilde{{\cal U} \Psi}^{1,k,s}$.
Moreover, as we have seen in section 6 after 6.4 that  
$P'(0) \in Comp_0 ~{\cal U} \Psi^{0,k,s}$, 
so $c'(0) \in Comp_0 ~\widetilde{{\cal U} \Psi}^{1,k,s}$.
It is easy to see that the map $c \mapsto c'(0)$ is
surjective and that the Lie bracket is the commutator bracket.
\hfill $\Box$

\me
Hence we have

\me 
\no
{\bf Theorem 7.14} {\it
The exact sequence of ILH Lie groups
\[ I \longrightarrow ({\cal U} \Psi^{0,k,s})_* \longrightarrow 
   ({\cal U} F^{0,k,s})_* \longrightarrow 
   {\cal D}^{s+k}_{\theta,0} (T^*M \setminus 0) 
   \longrightarrow e \]
has as corresponding sequence of ILH Lie algebras
\be
  0 \longrightarrow Comp_0~ {\cal U} \Psi^{0,k,s} 
  \longrightarrow  Comp_0~ \widetilde{{\cal U} \Psi}^{1,k,s} 
\stackrel{\varrho}{\longrightarrow}
  {\cal H}^1 \Omega^{0,2,s+k+1} (T^*M \setminus 0) 
  \longrightarrow 0 ,
\ee
where $\varrho = \frac{1}{i}\times$ principal symbol. }
\hfill $\Box$

\no
{\bf Remarks.} 

\no
1. As exact sequence of vector spaces (7.23) corresponds to 
\[ 0 \longrightarrow \sum\limits^k_{i=0} \Omega^{0,2,i+s} (S)   
   \longrightarrow \left( \sum\limits^k_{i=0} \Omega^{0,2,i+s} (S) \right)
   \oplus \Omega^{0,2,k+s+1} (S) \longrightarrow 
   \Omega^{0,2,k+s+1} (S) \longrightarrow 0 .
\]                                                   

\no
2. The space 
${\cal M}^\infty (I,B_\infty) = 
\lim\limits_{\begin{array}{c} \leftarrow \\[-2ex] r \end{array}} 
{\cal M}^r (I,B_\infty)$
splits into an uncountable number of components.
\[ {\cal M}^\infty (I,B_\infty) = \sum\limits_{i \in I} comp(g_i) , \]
where
\bea 
   comp(g) & = & \Big\{ g' | g' \en \mbox{satisfies} \en (I) \en \mbox{and} 
   \en (B_\infty), g \en \mbox{and} \en g' \en \mbox{are quasiisometric
   and} \en |g-g'|_{g,r} \nonumber \\
   & := & \Big( \int ( |g-g'|^2_{g,x} +
   \sum\limits^{r-1}_{i=0} |(\nabla^g)^i (\nabla^g-\nabla^{g'})
   |^2_{g,x}) dvol_x(g)  \Big)^{\frac{1}{2}} 
   < \infty \en \mbox{for all} \en r \Big\} . \nonumber
\eea
$g' \in comp(g)$ implies that 
$(S_g (T^* M), g_S)$ and $(S_{g'} (T^* M), g'_S)$,
are quasiisometric and $g'_S \in comp(g_S)$.
Moreover, all functional spaces which we considered and which enter
into the construction of
${\cal U} \Psi^{0,k,s}$, ${\cal U} F^{0,k,s}$, 
${\cal D}^sS_{\theta,0}$ are invariants of $comp(g)$ like the initial
spaces ${\cal U} \Psi^q$, ${\cal U} F^q$.
We obtain that the ILH Lie group structure of ${\cal U} F^{0,k}$ is
an invariant of $comp(g)$. We refer to \cite{Ei1} where we constructed the
spaces ${\cal M}^r(I,B_k)$.
\hfill $\Box$

\bi

\parbox[t]{3in}{Juergen Eichhorn\\                       
Institut fuer Mathematik\\
Universitaet Greifswald\\
Jahnstr. 15a\\
D-17487 Greifswald\\
Germany\\
eichhorn@ rz.uni-greifswald.de}      
\parbox[t]{2.5in}{Rudolf Schmid\\
Department of Mathematics\\
Emory University\\
Atlanta, Georgia 30322\\
USA\\
rudolf@ mathcs.emory.edu}                        
                         
\end{document}